\def\no{\if01}
\def\iftwelvept{\no}

\def\ifusepdf{\no}
\def\ifpsfont{\no}

\iftwelvept
\documentclass[12pt]{amsart}
\else
\documentclass[10pt]{amsart}
\fi
\usepackage{amssymb}
\usepackage{amscd}
\usepackage{latexsym}
\usepackage{verbatim}
\usepackage{mathrsfs}
\usepackage[all]{xy}

\numberwithin{equation}{section}

\ifusepdf
\usepackage{hyperref}
\else\fi
\ifpsfont
\usepackage[T1]{fontenc}
\usepackage{times}
\else\fi

\iftwelvept
\else\fi

\theoremstyle{plain}
\newtheorem{Theorem}{Theorem}[section]

\newtheorem{Proposition}[Theorem]{Proposition}
\newtheorem{Lemma}[Theorem]{Lemma}
\newtheorem{Corollary}[Theorem]{Corollary}

\theoremstyle{definition}

\newtheorem{Definition}[Theorem]{Definition}
\newtheorem{Remark}[Theorem]{Remark}
\newtheorem{Construction}[Theorem]{Construction}
\newtheorem{Example}[Theorem]{Example}

\newcommand{\ZZ}{\mathbf{Z}}

\newcommand{\DD}{\mathbb{D}}

\newcommand{\GG}{\mathcal{G}}
\newcommand{\HH}{\operatorname{\mathcal{HH}}}

\newcommand{\FF}{\mathcal{F}}

\newcommand{\BBB}{\mathcal{B}}

\newcommand{\DDD}{\mathcal{D}}

\newcommand{\CCC}{\mathcal{C}}

\newcommand{\Rep}{\operatorname{Rep}}

\newcommand{\Spec}{\operatorname{Spec}}

\newcommand{\Perf}{\operatorname{Perf}}

\newcommand{\SP}{\operatorname{Sp}}

\newcommand{\Mod}{\operatorname{Mod}}

\newcommand{\SSS}{\mathcal{S}}

\newcommand{\Map}{\operatorname{Map}}

\newcommand{\Fun}{\operatorname{Fun}}
\newcommand{\Alg}{\operatorname{Alg}}

\newcommand{\End}{\operatorname{End}}

\newcommand{\Aff}{\operatorname{Aff}}

\newcommand{\wCat}{\widehat{\textup{Cat}}_{\infty}}

\newcommand{\CAlg}{\operatorname{CAlg}}
\newcommand{\CCAlg}{\operatorname{CAlg}^{\le0}}
\newcommand{\CCAlgft}{\operatorname{CAlg}^{\le0,\diamondsuit}}
\newcommand{\CCAlgftec}{\operatorname{CAlg}^{\le0,\square}}

\newcommand{\QC}{\operatorname{QC}}

\newcommand{\EXT}{\operatorname{Art}^{\textup{tsz}}}
\newcommand{\TSZ}{\mathsf{TSZ}}

\newcommand{\FST}{\widehat{\mathsf{St}}^{\ast}}
\newcommand{\GFST}{\widehat{\mathsf{St}}^{!}}

\newcommand{\Ind}{\operatorname{Ind}}

\newcommand{\Coh}{\operatorname{Coh}}

\newcommand{\Free}{\operatorname{Free}}

\newcommand{\wSSS}{\widehat{\mathcal{S}}}

\newcommand{\eone}{\mathbf{E}_1}
\newcommand{\etwo}{\mathbf{E}_2}
\newcommand{\eenu}{\mathbf{E}_n}
\newcommand{\einf}{\mathbf{E}_\infty}

\newcommand{\KS}{\mathbf{KS}}
\newcommand{\Lie}{\mathbf{Lie}}

\newcommand{\LMod}{\operatorname{LMod}}
\newcommand{\RMod}{\operatorname{RMod}}
\newcommand{\ST}{\operatorname{\mathcal{S}t}}

\newcommand{\DR}{\textup{dR}}

\newcommand{\Def}{\operatorname{Def}}

\newcommand{\Proof}{{\sl Proof.}\quad}
\newcommand{\QED}{{\unskip\nobreak\hfil\penalty50\quad\null\nobreak\hfil
{$\Box$}\parfillskip0pt\finalhyphendemerits0\par\medskip}}

\begin{document}

\title[On D-modules of categories II]{On D-modules of categories II}

\author{Isamu Iwanari}

\maketitle

\section{Introduction}

This paper is the second in a series of papers.
Let $A$ be a smooth commutative (ordinary) algebra over a field of characteristic zero $k$.
Let $\CCC$ be a small $A$-linear stable $\infty$-category,
and let $\HH_\bullet(\CCC/A)$ be the Hochschild homology of $\CCC$ over $A$,
which is defined as an $A$-module spectrum
(equivalently, a differential graded (dg) $A$-module) endowed with an action of the circle
$S^1$.
In the first paper \cite{DC1}, 
we construct a lift of $\HH_\bullet(\CCC/A)$ to an $A\otimes_kS^1$-module spectrum
endowed with an $S^1$-action which is compatible with the $S^1$-action on $A\otimes_kS^1$.
Here $A\otimes_kS^1\simeq A\otimes_{A\otimes_kA}A$.
Since $\Spec A\otimes_kS^1$ may be regarded as the (derive) loop space of $S=\Spec A$,
such a lift can be thought of as a deformation/extension of $\HH_\bullet(\CCC/A)$
along $S=\Spec A\to LS=\Spec A\otimes_kS^1$ with respect to $S^1$-actions.
Using a lift to $LS$ we construct a $D$-module structure on the periodic cyclic homology/complex
$\mathcal{HP}_\bullet(\CCC/A)$.
We provided two methods for constructing a lift.

(I) The first method uses the canonical extension of factorization homology to mapping stacks.
When $\CCC$ is an stable idempotent-complete $\infty$-catergory
over a scheme $S$,
there exists a canonically defined $S^1$-equivariant extention $\HH_\bullet(\CCC/S)_L$ of the relative Hochschild homology (chain complex)
$\HH_\bullet(\CCC/S)$ to the (derived) loop space $LS$.
This approach is simple and easy, and the resulting object has a nice functoriality.
Moreover, it admits a vast generalization to the lifts of factorization homology of $\eenu$-algebras.

(II) The second method uses the algebra of the pair $(\HH^\bullet(\CCC/A),\HH_\bullet(\CCC/A))$
of Hochschild cohomology and Hochschild homology (by which we mean a chain complex/spectra
computing Hochschild cohomology and Hochschild homology).
The algebraic structure may be thought of as a version of Cartan calculus and is defined as an algebra over a colored topological operad called Kontsevich-Soibelman operad.
An advantage of this approach is a direct relation with the pair $(\HH^\bullet(\CCC/A),\HH_\bullet(\CCC/A))$. This relation is useful.
For example, if we write $\mathcal{H}_{\circlearrowleft}(\CCC)\in \QC_!(LS)^{S^1}$ for the lift, then
the pullback of $\mathcal{H}_{\circlearrowleft}(\CCC)$ to $\QC_!(S\times_kLS)^{\wedge}_S)$ can
be described in terms of the (dg) Lie module coming from the Lie derivation and the contraction map
built in the pair $(\HH^\bullet(\CCC/A),\HH_\bullet(\CCC/A))$ together with the Kodaira-Spencer morphism
for $\CCC$. Here $(S\times_kLS)^{\wedge}_S$ is the formal stack obtained from $S\times_kLS$ by the formal completion (see Section~\ref{formalstacksec}).
This structure will be applied to the study of the resulting object in the future
work.

Both have their own pleasant features
so that it is desirable to compare them.
The main result of this paper is a comparison of the two methods.
We state the main result in a naive way (see Theorem~\ref{main} for the precise statement): 

\begin{Theorem}
Assume that $S$ is affine and smooth over $k$.
The two lifts constructed by the two methods (I) and (II) coincide. Namely,
$\mathcal{H}_\circlearrowleft(\CCC)$ 
and $\HH_\bullet(\CCC/S)_{L}=\HH_\bullet(\CCC/k)$ coincide in 
the $\infty$-category $\QC_!(LS)^{S^1}$ of $S^1$-equivariant 
Ind-coherent complexes on the derived loop space $LS$.
\end{Theorem}

The remarkable feature of this comparison result is that it can be thought of as Koszul duality
between two methods.
The first method (I) mainly uses module objects over commutative algebras ($\einf$-algebras).
On the other hand, in the second method (II), the main data is dg Lie algebra modules over dg Lie algebras.
To compare two methods, it is necessary to relate modules appearing in (I) and dg Lie algebra modules in (II).
One of the key relations is Koszul duality between modules over an argmented commutative algebra $B$
and modules over Koszul dual dg Lie algebra $\DD_\infty(B)$: results and machinery
in \cite{IMA} play an important role.

\section{Preliminaries}

Throughout this paper, $k$ is a field of characteristic
zero and $A$ is a commutative smooth (ordinary) algebra
over $k$.

\subsection{Convention and Notation}
\label{NC}
We use the theory of $(\infty,1)$-categories.
The principal model of $(\infty,1)$-categories in this paper
is the theory of quasi-categories entensively developed in \cite{HTT}, \cite{HA}.
Following \cite{HTT}, we call quasi-categories $\infty$-categories.
We use the notation in \cite{DC1}.
In particular, we frequently use the following symbols:

\begin{itemize}

\item $\Delta^n$: the standard $n$-simplex


\item $\SSS$: $\infty$-category of small spaces/$\infty$-groupoids. We denote by $\widehat{\SSS}$
the $\infty$-category of spaces in the enlarged universe (cf. \cite[1.2.16]{HTT}).

\item $\CCC^\simeq$: the largest Kan subcomplex of an $\infty$-category $\CCC$. Namely, $\CCC^\simeq$ is the largest $\infty$-groupoid contained in $\CCC$.

\item $\CCC^{op}$: the opposite $\infty$-category of an $\infty$-category. We also use the superscript ``op" to indicate the opposite category for ordinary categories and enriched categories.

\item $\Map_{\mathcal{C}}(C,C')$: the mapping space from an object $C\in\mathcal{C}$ to $C'\in \mathcal{C}$ where $\mathcal{C}$ is an $\infty$-category.
We usually view it as an object in $\mathcal{S}$ (cf. \cite[1.2.2]{HTT}).

\item $\operatorname{Cat}_\infty$: the $\infty$-category of small $\infty$-categories.
$\wCat$ is the $\infty$-cartegory of large $\infty$-categories.


\item $\Fun(A,B)$: the function complex for simplicial sets $A$ and $B$. If $A$ and $B$ are $\infty$-categories, we regard $\Fun(A,B)$ as the functor category. 

\item $\CAlg_R$: the $\infty$-category of $R$-module spectra where $R$ is a commutative
algebra object (i.e. an $\eenu$-algebra object) in the $\infty$-category $\SP$ of spectra.
For a symmetric monoidal $\infty$-category $\mathcal{M}$ we write $\CAlg(\mathcal{M})$
for the $\infty$-category of commutative algebra objects in $\mathcal{M}$.

\item $\Mod_B$: the (symmetric monoidal) $\infty$-category of $B$-module obejcts in $\SP$
for $B\in \CAlg(\SP)$.
We also write $\QC(B)$ or $\QC(\Spec B)$ for $\Mod_B$. Namely, $\QC(B)=\QC(\Spec B)=\Mod_B$.
We denote by $\Perf_R$ the full subcategory of $\Mod_B$, which consists of 
dualizable objects.
For $B\in \CAlg_R=\CAlg(\Mod_R)$,
the forgetful functor induces an equivalence $\Mod_B(\Mod_R)\simeq \Mod_B(\SP)$.
If $R\in \CAlg_k$ is connective, we write $\CCAlg_R$ for the full subcategory
of $\CAlg_R$ spanned by connective objects with respect to the standart $t$-structure
on $\Mod_R$. 
Let $\CCAlgft_k$  denote the full subcategory of $\CAlg_k$ spanned by
connective objects almost of finite type over $k$.
Let $\CCAlgftec_k$ denote the full subcategory of $\CCAlgft_k$ spanned by
those objects $R$ such that $H_n(R)=0$ for $n>>0$.

\item $\QC:\CAlg_k\to \wCat$ : the functor which carries $B\in \CAlg_k$ to
$\QC(B)$.
This functor sends $B\to B'$ (corresponding to $f:\Spec B'\to \Spec B$)
to the $*$-pullback functor $f^*:\QC(\Spec B')\to \QC(\Spec B)$.

\item $\QC_!:\CCAlgft_k\to \wCat$: the functor which carries $B\in \CCAlgft_k$ to
the $\infty$-category $\QC_!(B)=\Ind(\Coh(B))$ of Ind-coherent sheaves/complexes over $\Spec B$.
This functor sends $B\to B'$ to the $!$-pullback functor $f^!:\QC_!(\Spec B')\to \QC_!(\Spec B)$.
See \cite{Ind}.

\item $\Upsilon:\QC|_{\CCAlgft_k}\to \QC_!$: the natural transformation
which induces $\Upsilon_B:\QC(B)\to \QC_!(B)$ given by $M\mapsto M\otimes_B\omega_B$
for each $B$, where $\omega_B$ is the $!$-pullback $p^!(k)$ of $k\in \QC_!(k)$ along
$p:\Spec B\to \Spec k$.
See \cite{Ind}.

\end{itemize}

\subsection{Formal stacks}

We put $S=\Spec A\in \Aff_k=(\CCAlg_k)^{op}$.

\label{formalstacksec}

As in \cite{DC1} we use the theory of pointed formal stacks over $A$,
which was developed in \cite{H}, \cite[Vol.II]{Gai2}. The theory generalizes the theory of formal
moduli problems developed in \cite[X]{DAG}, which should be thought of as the theory
of pointed formal stacks over $k$.
We refer the reader to \cite{H}, \cite{Gai2}, \cite{DAG} for detail and to \cite[Section 3.4]{DC1} for 
a brief review.

Let $Lie_A$ be the $\infty$-category of dg Lie algebras.
The $\infty$-category $Lie_A$ is obtained from the model category of dg Lie algebras 
(whose fibrations are termwise surjective maps) by inverting
quasi-isomorphisms
(another equivalent approach is to define it as
the $\infty$-category obtained from algebras over the Lie operad $\Lie$).
Let
$\EXT_A$ be the full subcategory of $(\CCAlg_k)_{A//A}:=((\CCAlg_k)_{A/})_{/A}$,
which is spanned by  trivial square zero extensions $A=A\oplus 0\hookrightarrow A\oplus M\stackrel{p_1}{\to} A$ such that
$M$ is a connective $A$-module of the form $\oplus_{1\le i\le n}A^{\oplus r_i}[d_i]$ ($r_i\ge0,\ d_i\ge0)$.
We note that any object $R\to A$ of $\EXT_{A}$ is a trivial square zero extension of the form $\textup{pr}_1:A\oplus M \to A$
such that $M=\oplus_{1\le i\le n}A^{\oplus r_i}[p_i]$  ($p_i\ge0$).
Thus $\Map_{\CAlg_k}(A\oplus M, A)$ is a contractible space
since $A$ is an ordinary commutative reduced algebra over $k$.
It follows that the composite 
functor
$\EXT_A\to (\CCAlg_k)_{A//A}\to (\CCAlg_k)_{A/}\simeq \CCAlg_A$
is fully faithful. 
We can also think of $\EXT_A$ as a full subcategory of $\CCAlg_A$.

By abuse of notation, we often write $R$ for an object $A\to R\to A$ of $\CCAlg_{A//A}$.
Similarly, we often omit the augmentations from the notation.
Let $\TSZ_A$ denote the opposite category of $\EXT_A$.
A pointed formal stack over $A$ is a functor $\EXT_A\to \SSS$
 satisfying a certain ``gluing condition'' (cf. \cite{H}, \cite[Section 3.4]{DC1}).
The $\infty$-category $\FST_A$ of  is
 a full subcategory of $\Fun(\EXT_A,\SSS)$.
The Yoneda embedding $\TSZ_A\hookrightarrow \Fun(\EXT_A,\SSS)$ factors through
$\FST_A\subset \Fun(\EXT_A,\SSS)$, and we often regard $\TSZ_A$ as a full subcategory of $\FST_A$.
Let $\Free_{Lie}:\Mod_A\to Lie_A$ be the free Lie algebra functor
which is a left adjoint to the forgetful functor $Lie_A\to \Mod_A$. 
Let $\Mod_{A}^{f}\subset \Mod_A$ be the full subcategory that consists of
objects of the form $\oplus_{1\le i\le n}A^{\oplus r_i}[d_i]$ ($d_i\le -1$)
Let $Lie_A^f$ be the full subcategoy of $Lie_A$, which is the essential image of the restriction of
the free Lie algebra functor  $\Mod_A^{f}\to Lie_A$.
Thanks to \cite[1.5.6]{H}, this adjoint pair induces an adjoint pair
\[
\xymatrix{
\mathcal{F}:Lie_A  \ar@<0.5ex>[r]^{\simeq} &   \FST_A :\mathcal{L}.\ar@<0.5ex>[l]^{\simeq}  
}
\]
which are inverse to one another.
This adjoint pair extends the pair of mutually inverse functors $Ch^\bullet:Lie_A^f\simeq \TSZ_A:\DD_\infty$,
which is the restriction of the Koszul duality adjoint pair
$Ch^{\bullet}:Lie_A\rightleftarrows (\CAlg_{A//A})^{op}:\DD_{\infty}$ such that
the left adjoint $Ch^\bullet$ is the Chevalley-Eilenberg cochain functor.
For $L\in Lie_A$, we usually write $\mathcal{F}_L$ for the associated formal stack $\mathcal{F}(L)\in \FST_A$.

\begin{Example}
Let $\widehat{S\times_kS}$ be the pointed formal stack obtained from
the formal completion of $S\times_kS\stackrel{\textup{pr}_1}{\to} S$
along the diagonal $\Delta:S\to S\times_kS$.
The pointed formal stack $\widehat{S\times_kS}$
is defined to be the functor
$\EXT_A\to \SSS$
given by $[A\to R\to A]\mapsto \Map_{(\Aff_k)_{S//S}}(\Spec R,S\times_kS)\simeq \Map_{(\Aff_k)_{S/}}(\Spec R, S)$.
Let $LS$ denote the (derived) loop space defined by $LS=\Spec A\otimes_kS^1$.
It has the obvious $S^1$-action and the $S^1$-equivariant morphism
$\iota:S\to LS$ determined by $A\otimes_kS^1\to A\otimes_k\ast\simeq A$ induced by
the contraction $S^1\to \ast$, where the $S^1$-action on $\ast$ is trivial.
Let $\widehat{S\times_kLS}$ be the pointed formal stack obtained from
$S\times_kLS\stackrel{\textup{pr}_1}{\to} S$
along $\textup{id}\times \iota:S\to S\times_kLS$ by the formal completion defined in the same way.
\end{Example}

Let $\GFST_A$ denote the $\infty$-category of pointed formal moduli problems
over $A$ defined in \cite{Gai2}. 
It can be considered as a full subcategory of $(\Fun(\CCAlgftec_k,\SSS))_{\Spec A//\Spec A}$.
We will dub an object of $\GFST_A$ as a pointed $!$-formal stack over $A$ (or $S$).
The Yoneda embedding $\TSZ_A\to (\Fun(\CCAlgftec_k,\SSS))_{\Spec A//\Spec A}$
factors through $\GFST_A$.
There exists a categorical equivalence $\GFST_A\simeq Lie_A$.
See \cite[Vol. II, Chap.5]{Gai2} (also \cite[Section 7.1]{DC1} for the quick review).
According to \cite[Construction 7.2, Remark 7.4]{DC1},
there eixsts a categorical equivalence $\Theta_A:\FST_A\stackrel{\sim}{\to}\GFST_A$
which commutes with $\TSZ_A\hookrightarrow \FST_A$ and $\TSZ_A\hookrightarrow \GFST_A$.
We remark that this equivalence is due to the smoothness of $A$.

\begin{Example}
Let $(S\times_kS)_{S}^\wedge$ denote $(S\times_kS)\times_{(S\times_kS)_{\DR}}S_{\DR}$,
that is determined by the diagonal $\Delta:S\to S\times_k S$
(see \cite[Vol.II, Chap.4]{Gai2} or \cite[the review after Remark 7.4]{DC1}).
By definition, $(S\times_kS)_{S}^\wedge$ is the functor
$(\CCAlgftec_k)_{A//A} \to \SSS$ defined by $R\mapsto (S(R)\times S(R))\times_{(S(R_{red})\times S(R_{red}))}S(R_{red})$. 
We think of it as
the pointed $!$-formal stack 
obtained from $S\times_kS$ by taking the formal completion along the diagonal $S\to S\times_kS$.
Let $LS=\Spec A\otimes_kS^1=S\times_{S\times_kS}S$ be the free loop space of derived scheme $S=\Spec A$ over $k$.
Let $(S\times_kLS)_{S}^\wedge$ denote $(S\times_kLS)\times_{(S\times_kLS)_{\DR}}S_{\DR}$,
determined by $\textup{id}\times \iota:S\to S\times_k LS$.
The equivalence $\Theta_A:\FST_A \stackrel{\sim}{\to} \GFST_A$ carries $\widehat{S\times_kS}$ and $\widehat{S\times_kLS}$
to $(S\times_kS)^\wedge_{S}$ and $(S\times_kLS)^\wedge_{S}$, respectively (see \cite[Proposition 7.7]{DC1}). Thus, we will regard $(S\times_kS)^{\wedge}_S$ and $(S\times_k LS)^{\wedge}_S$
as the images of $\widehat{S\times_kS}$ and $\widehat{S\times_kLS}$, respectively.
\end{Example}

\section{The absolute Hochschild homology and cyclic deformations}

\subsection{}
Let $\CCC$ be a small stable idempotent-complete $\infty$-category.
Let $\ST$ denote the $\infty$-category which consists of small stable 
idempotent-complete $\infty$-categories where mapping spaces are 
spanned by exact functors (see e.g. \cite{BGT1}, \cite[Section 2]{DC1} for detail).
There exists a closed symmetric monoidal structure on $\ST$
such that the internal Hom/mapping object is given by $\Fun^{\textup{ex}}(-,-)$.
Here  $\Fun^{\textup{ex}}(-,-)$ indicates the full subcategory of $\Fun(-,-)$
spanned by exact functors.
For $R\in \CAlg(\SP)$
we consider the symmetric monoidal stable $\infty$-category $\Perf_R$ as a commutative algebra
object in the symmetric monoidal $\infty$-category $\ST$.
We define $\ST_R$ to be $\Mod_{\Perf_R}(\ST)$ and refer to an object of $\ST_R$
as an $A$-linear small stable (idempotent-complete) $\infty$-category.
Namely, by an $A$-linear structure on $\CCC$ we mean a $\Perf_R$-module structure on $\CCC$.
Moreover, $\ST_R=\Mod_{\Perf_R}(\ST)$ inherits a symmetric monoidal structure
from that on $\ST$.

\subsection{}

\label{cocartesianfamily}

Let $\HH_\bullet(-/A):\ST_A^\otimes\to (\Mod_A^{S^1})^{\otimes}$ be the symmetric monoidal functor which carries an $A$-linear stable small stable $\infty$-category
$\DDD$ to the Hochschild homology $A$-module spectrum $\HH_\bullet(\DDD/A)$ (see \cite[Section 6]{I}).
We define the symmetric monoidal functor $\HH_\bullet(-/k):\ST_k^\otimes\to \Mod_k^\otimes$
in the same way.

We apply the construction
$\Mod(\mathcal{M})\to \CAlg(\mathcal{M})$ (see \cite[Chaper 3, Theorem 4.5.3.1]{HA}) to 
$\HH_\bullet(-/A):\ST_A\to \Mod_A^{S^1}$ to obtain
\[
\xymatrix{
\Mod(\ST_A) \ar[rr]^(0.45){\Mod(\HH_\bullet(-/A))} \ar[d] & & \Mod(\Mod_A^{S^1}) \ar[d] \\
\CAlg(\ST_A) \ar[rr]^(0.45){\CAlg(\HH_\bullet(-/A))}  & & \CAlg(\Mod_A^{S^1})
}
\]
where vertical functors coCartesian fibrations.

Let $\CCC$ be an object of $\ST_A$.
We define $\Mod^{\textup{alg}}(\ST_A)$ to be the fiber product
$\Mod(\ST_A)\times_{\CAlg(\ST_A)}\CAlg_A$, which is 
determined by the base change along
$\CAlg_A=\CAlg(\Mod_A) \to \CAlg(\ST_A^\otimes)$ which carries
$R$ to $\Perf_R^\otimes$ (it is obtained from the construcion in \cite[4.8.5.21]{HA}, see also \cite[Section 2]{DC1}). 
We set $\Mod^{\textup{alg}}(\ST_A)^+:=\Mod^{\textup{alg}}(\ST_A)\times_{\CAlg_A}\CAlg_A^+$.
Put $\Mod^{\textup{alg}}(\ST_A)_{\CCC}:=\Mod^{\textup{alg}}(\ST_A)^+\times_{\ST_A}\{\CCC\}$.
Let $\Mod^{\textup{alg}}(\ST_A)\times_{\CAlg_A}\CAlg_A^+ \to \ST_A$ be the functor
which carries $\CCC_R \in \Mod_{\Perf_R}(\ST_A)$
to $\CCC_R\otimes_{\Perf_R}\Perf_A$ (see \cite[Section 5.1]{IMA} for the construction
up to $\Mod(\ST_A)\simeq \RMod(\ST_A)\times_{\Alg_1(\ST_A)}\CAlg(\ST_A)$).
 Consider the coCartesian fibration
\[
\Mod^{\textup{alg}}(\ST_A)_{\CCC}:=\{\CCC\}\times_{\ST_A}\Mod^{\textup{alg}}(\ST_A)\times_{\CAlg_A}\CAlg_A^+\to \CAlg_A^+.
\]
This coCartesian fibration
corresponds to the functor $\Def_{\einf}'(\CCC) : \CAlg_A^+\to \wCat$
 informally given by $[R\to A]\mapsto \ST_R\times_{\ST_A}\{\CCC\}$.
If $g:\wCat\to \wSSS$ denotes the functor defined by $g(V)=V^{\simeq}$
(that is obtained by taking the largest $\infty$-groupoid/Kan complex contained in $V$),
the composition with $g$ determines $\Def_{\einf}(\CCC)=g\circ \Def_{\einf}'(\CCC): \CAlg_A^+\to \wSSS$.
We write $\Def_{\CCC}:\EXT_A \to \wSSS$ for the composite
$\EXT_A\hookrightarrow \CAlg_A^+ \stackrel{\Def_{\einf}(\CCC)}{\longrightarrow} \wSSS$.
Let $X:\EXT_A\to \SSS$ be a functor, that is, a pointed formal prestack.
 Since $X$ is a colimit of $(\TSZ_A)_{/X}\to \TSZ_A\to \Fun(\EXT_A,\SSS)$, there exists
 a canonical equivalence
\[
\Map_{\Fun(\EXT_A,\wSSS)}(X,\Def_{\CCC})\simeq \lim_{\Spec R\in (\TSZ_A)_{/X}}\Def_{\CCC}(R).
\]

\subsection{}
\label{cyclicdeformationsec}
Let $\mathcal{H}$ be an $A$-module spectrum endowed with an $S^1$-action,
that is, an object of $\Mod_A^{S^1}=\Fun(BS^1,\Mod_A)$.
We briefly review the notion of cyclic deformations of $\mathcal{H}$.
Set the coCartesian fibration
\[
\Mod(\Mod_A^{S^1})_{\mathcal{H}}:=\{\mathcal{H}\}\times_{\Mod_A^{S^1}}\Mod(\Mod_A^{S^1})\times_{\CAlg_A^{S^1}}(\CAlg_A^{S^1})_{/A}\to (\CAlg_A^{S^1})_{/A}.
\]
Here $\CAlg_A^{S^1}\simeq \Fun(BS^1, \CAlg_A)\simeq \CAlg(\Mod_A^{S^1})$,
and $(\CAlg_A^{S^1})_{/A}$ denotes the overcategory over  the unit algebra $A$ endowed
with the trivial $S^1$-action.
This corresponds to the functor $(\CAlg_A^{S^1})_{/A}\to \wCat$
given by $[B\to A]\mapsto \Mod_{B}(\Mod_A^{S^1})\times_{\Mod_A^{S^1}}\{\CCC\}$.
Let $\EXT_A\to (\CAlg_A^{S^1})_{/A}$ is the composite
$\EXT_A\stackrel{\textup{forget}}{\longrightarrow} (\CAlg_A)_{/A}\stackrel{\otimes_AS^1}{\longrightarrow} (\CAlg_A^{S^1})_{/A}$
(which is nothing but $\EXT_A\stackrel{T}{\to} \Fun(BS^1,(\CCAlg_k)_{A//A})\subset (\CAlg_A^{S^1})_{/A}$ in Section~\ref{completion}).
Let $\Def^{'\circlearrowleft}(\mathcal{H}):\EXT_A \to \wCat$ be the functor
corresponding to the base change
$\Mod(\Mod_A^{S^1})_{\mathcal{H}}\times_{ (\CAlg_A^{S^1})_{/A}}\EXT_A\to \EXT_A$.
As in the case of $\Def_{\CCC}$, we define $\Def^{\circlearrowleft}(\mathcal{H}):\EXT_A \to \wSSS$
to be $g\circ \Def^{'\circlearrowleft}(\mathcal{H}):\EXT_A \to \wCat\to \wSSS$.
We refer to $\Def^{\circlearrowleft}(\mathcal{H})(R)$ as the $\infty$-groupoid/space
of cyclic deformations of $\mathcal{H}$ to $R$. 
Suppose that $X:\EXT_A\to \SSS$ is a functor. There exists a canonical equivalence
\[
\Map_{\Fun(\EXT_A,\wSSS)}(X,\Def^\circlearrowleft(\mathcal{H}))\simeq \lim_{\Spec R\in (\TSZ_A)_{/X}}\Mod_{R\otimes_AS^1}(\Mod_A^{S^1})^\simeq\times_{(\Mod_A^{S^1})^\simeq}\{\mathcal{H}\}.
\]

\begin{Definition}
\label{QCHcircle}
Let $\QC:\CAlg_A\to \wCat$ for the functor
corresponding to the coCartesian fibration $\Mod(\Mod_A)\to \CAlg_A$ (cf. Section~\ref{NC}, \cite[4.5.3.1]{HA}).
Then 
we define $\QC_H^\circlearrowleft$ to be the composite functor
\[
\xymatrix{
\EXT_A\ar[r]^(0.5){\textup{forget}} & \CAlg_A \ar[r]^(0.4){\otimes_AS^1}  &   \Fun(BS^1,\CAlg_A)  \ar[r]^{\QC} & \Fun(BS^1,\wCat).
}
\]
We also define $\QC_H^\circlearrowleft:\Fun(\EXT_A,\SSS)^{op}
\to \Fun(BS^1,\wCat)$
to be the right Kan extension of $\QC_H:\EXT_A\to \Fun(BS^1,\wCat)$
along the Yoneda embedding $\EXT_A\to \Fun(\EXT_A,\SSS)^{op}$
(we abuse notation by using the same symbol $\QC_H^\circlearrowleft$).
Using this functor we regard $\Def^\circlearrowleft(\mathcal{H})(R)$
as $\bigl(\QC_H^\circlearrowleft(R)^{S^1}\times_{\QC_H^{\circlearrowleft}(A)^{S^1}}\{\mathcal{H}\}\bigr)^{\simeq}$.
\end{Definition}

\subsection{}
Take $\mathcal{H}$ to be $\HH_\bullet(\CCC/A)$.
The symmetric monoidal functor
$\HH_\bullet(-/A):\ST^\otimes_A\to (\Mod_A^{S^1})^\otimes$
gives rise to the diagram
\[
\label{diagramDD}
\vcenter{
\xymatrix{
\Mod^{\textup{alg}}(\ST_A)_{\CCC} \ar[r] \ar[d] &  \Mod(\Mod_A^{S^1})_{\mathcal{H}} \ar[d] \\
\CAlg_A^{+} \ar[r] & (\CAlg_A^{S^1})_{/A}
}\tag{3.1}}
\]
By \cite[Lemma 6.3]{IMA}, $\Mod(\ST_A)\to \Mod(\Mod_A^{S^1})$
preserves coCartesian morphisms. It follows that 
$\Mod^{\textup{alg}}(\ST_A)_{\CCC} \to  \Mod(\Mod_A^{S^1})_{\mathcal{H}}$
 preserves coCartesian morphisms.

Write $h:=\HH_\bullet(-/A):\CAlg_A\simeq \CAlg(\Alg_1(\Mod_A)) \to \CAlg(\Mod_A^{S^1})\simeq \CAlg_A^{S^1}$.
Since $h:\CAlg_A\to \CAlg_A^{S^1}$ is equivalent to 
the functor $\otimes_AS^1:\CAlg_A\to \CAlg_A^{S^1}$ given by the tensor by $S^1$
(see \cite[Lemma 3.5]{IMA}), 
it follows that if $i:\CAlg_A\to \CAlg_A^{S^1}$ is the functor
which sends each $B$ to $B$ with the trivial $S^1$-action there exists
the natural transformation $h\to i$ induced by the contraction $S^1\to \ast$ to the one-point
space.
For each $B\in \CAlg_A$, $h(B)\to i(B)$
can be identified with $B\otimes_AS^1\to B\otimes_A\ast=B$.
The natural transformation $h\to i$ naturally extends to
a natural transformation $h^+\to i^+$ between functors
$\CAlg_A^+\simeq \CAlg(\Alg_1(\Mod_A)) \to \CAlg(\Mod_A^{S^1})_{/A}$
where $h^+$ and $i^+$ are functors induced by $h$ and $i$ in the obvious way.

We write $U_A$ for the full subcategory of $\bigl((\Aff_k\bigr)_{S//S})_{/S\times_kS}$
which consists of those objects $\Spec C\to S\times_kS$ such that $\Spec C\in \TSZ_A$.
Here $S\times_kS$ in $(\Aff_k)_{S//S}$ indicates 
$S\stackrel{\Delta}{\to} S\times_kS\stackrel{\textup{pr}_1}{\to}S$
having the diagonal morphism followed by the first projection.  
Let $\overline{U}_A$ be the full subcategory of $\bigl((\Aff_k\bigr)_{S//S})_{/S\times_kS}$
which
is obtained from $U_A$ by adding the object $\textup{id}:S\times_kS\to S\times_kS$.
We note that there exists an equivalence $\overline{U}_A\simeq U_A^{\triangleright}$,
provided that $\dim A>0$.
Let $\overline{\rho}:\overline{U}_A^{op}\to  \CAlg_A^+$ be the forgetful functor.
The natural transformation $h^+\to i^+$ induces
$h^+\circ \overline{\rho}\to i^+\circ \overline{\rho}$ which is described as $\tau:\Delta^1\times \overline{U}_A^{op}\to \CAlg(\Mod_A^{S^1})_{/A}$.

\begin{Proposition}
\label{cocarsection}
Let $\textup{pr}_2^*(\CCC)$ denote the base change of $\CCC$ along
$\textup{pr}_2: S\times_kS \to S$, that is, $\CCC\otimes_{\Perf_A}\Perf_{(A\otimes_kA)}$.
Then there exists an essentially unique functor $\sigma$ filling the diagram
\[
\xymatrix{
 &  \Mod(\Mod_A^{S^1})_{\mathcal{H}} \ar[d]^\pi \\
\Delta^1 \times \overline{U}_A^{op} \ar[r]_(0.45){\tau} \ar[ru]^{\sigma}& (\CAlg_A^{S^1})_{/A}
}
\]
such that
\begin{enumerate}

\item $\sigma$ sends the initial object $(0,S\times_kS)$ of $\Delta^1 \times \overline{U}_A^{op}$ to
$\HH_\bullet(\textup{pr}_2^*(\CCC)/A)$,

\item the functor
$\sigma:\Delta^1\times \overline{U}_A^{op}\to  \Mod(\Mod_A^{S^1})_{\mathcal{H}}$
sends any morphism to a coCartesian morphism.
\end{enumerate}
\end{Proposition}

We start with the following Lemmata.

\begin{Lemma}
\label{coCartesiantransformation}
Let $\pi:\mathcal{P}\to \BBB$ be a coCartesian fibration between $\infty$-categories.
Let $f:I\to \BBB$ and $g:I\to \BBB$ be functors between $\infty$-categories.
Suppose that we are given a natural transformation $\sigma:\Delta^1\times I\to \BBB$
from $f$ to $g$.
Let $u:I\to \mathcal{P}$ be a functor such that $f=\pi\circ u$.
Then there exists an essentially unique $\overline{u}:\Delta^1\times I \to \mathcal{P}$
which extends $u:\{0\}\times I=I\to \mathcal{P}$ such that 
(i) $\sigma=\pi\circ \overline{u}$,
and (ii) 
for any object $x\in I$
the induced map $\Delta=\Delta^1\times\{x\} \to \mathcal{P}$ determines 
a coCartesian morphism/edge in $\mathcal{P}$.
\end{Lemma}

\Proof
Consider $\pi^I:\Fun(I,\mathcal{P})\to \Fun(I,\BBB)$
induced by $\sigma$.
According to \cite[3.1.2.1 (1)]{HTT}
this functor is a coCarrtesian fibration.
Passing to adjoints, we have
$v:\Delta^0\to \Fun(I,\mathcal{P})$
and $\tau:\Delta^1\to \Fun(I,\BBB)$ 
which correspond to $u$ and $\sigma$, respectively.
Then there exists an essentially unique coCatesian morphism $\overline{v}:\Delta^1\to\Fun(I, \mathcal{P})$
which lies over $\tau$ and extends $v$.
Then according to \cite[3.1.2.1 (2)]{HTT}
$\overline{u}:\Delta^1\times I\to \mathcal{P}$ corresponding to $\overline{v}$
satisfies the required property (ii).
\QED

\begin{Lemma}
\label{cocartesianextension}
Let $f: I\to \BBB$ be a functor between $\infty$-categories.
Let $\pi:\mathcal{P}\to \BBB$ be a coCartesian fibration.
Suppose that $s$ is an initial object of $I$.
Let $\Fun_{\BBB}(I, \mathcal{P})$ denote the function complex over $\BBB$.
Let $\Fun_{\BBB}^\dagger(I, \mathcal{P})$
be the full subcategory of $\Fun_{\BBB}(I, \mathcal{P})$
spanned by $u:I\to \mathcal{P}$ such that $u(\alpha)$ is a coCartesian 
morphism for any morphism $\alpha$ in $I$.
The evaluation at $s$ induces
an equivalence $\Fun_{\BBB}^\dagger(I, \mathcal{P}) \stackrel{\sim}{\to} \pi^{-1}(f(s))$.
\end{Lemma}

\Proof
Let $\theta:\BBB\to \wCat$ be the functor corresponding to the coCartesian fibration
$\pi$
through the straightening functor (cf. \cite[3.2]{HTT}).
According to \cite[3.3.3.2]{HTT}, $\Fun_{\BBB}^\dagger(I, \mathcal{P})$ can be identified with
a limit of $\theta\circ f:I\to \BBB\to \wCat$.
Since $s$ is an initial object,
$\{s\}^{op}\to I^{op}$ is cofinal.
It follows that  $\Fun_{\BBB}^\dagger(I, \mathcal{P}) \to \pi^{-1}(f(s))$
 is a categorical equivalence.
 \QED

{\it Proof of Proposition~\ref{cocarsection}.}
By Lemma~\ref{cocartesianextension},
there exists an essentially unique $\overline{\rho}':\overline{U}_A^{op}\to \Mod^{\textup{alg}}(\ST_A)_{\CCC}$ over $\CAlg_A^+$,
which carries the initial object to $\textup{pr}_2^*(\CCC)$ lying over $A\otimes_kA$
and carries
 any morphism in $\overline{U}_A^{op}$ to a coCartesian morphism.
The composition with horizontal functors in the diagram~\ref{diagramDD} gives rise to 
$\sigma_0:\{0\}\times \overline{U}_A^{op}\to  \Mod(\Mod_A^{S^1})_{\mathcal{H}}$
such that $\pi\circ \sigma_0=\tau|_{\{0\}\times \overline{U}_A^{op}}=h^+\circ \overline{\rho}$, and each morphism in $\{0\}\times\overline{U}_A^{op}$ maps to a coCartesian morphism
$\Mod(\Mod_A^{S^1})_{\mathcal{H}}$ and the image of $(0,S\times_kS)$ is $\HH_\bullet(\textup{pr}_2^*(\CCC)/A)$.

Applying Lemma~\ref{coCartesiantransformation} to the natural transformation $\tau$
and the diagram $\sigma_0$, we obtain $\sigma:\Delta^1\times \overline{U}_A^{op}\to  \Mod(\Mod_A^{S^1})_{\mathcal{H}}$ 
such that
$\sigma$ lies over $\tau$,
and
for each $T\in \overline{U}_A^{op}$, $\Delta^1\times \{T\}\to \Mod(\Mod_A^{S^1})_{\mathcal{H}}$ determines a coCartesian morphism.
The functor $\sigma$ is unique in the sense that given the fixed image $\HH_\bullet(\textup{pr}_2^*(\CCC)/A)$ of  $(0,S\times_kS)$,
the $\infty$-category classifying $\sigma$ having property (2) is the contractible space
(see Lemma~\ref{cocartesianextension}).
To see this, it is enough to prove that
the functor
$\sigma$
sends any morphism to a coCartesian morphism.
Let $f:\Delta^1\to \Delta^1\times \overline{U}_A^{op}$ be a morphism in $\Delta^1\times \overline{U}_A^{op}$.
If $f$ factors through $\{0\}\times \overline{U}_A^{op}$, then the assertion is obvious
since $\overline{\rho}'$ maps any morphism to a coCartesian morphism, and
$\Mod^{\textup{alg}}(\ST_A)_{\CCC} \to  \Mod(\Mod_A^{S^1})_{\mathcal{H}}$
preserves coCartesian morphisms.
Next, we consider the case when the the source $s$ of $f$ is contained in $\{0\}\times \overline{U}_A^{op}$
and the target $t$ is contained in $\{1\}\times \overline{U}_A^{op}$.
If $f_0$ denotes the morphism $s_0\to t_0$ in $\overline{U}_A^{op}$ determiend by the composite $\Delta^1\to\Delta^1\times \overline{U}_A^{op}\to \overline{U}_A$, then $f$ is equivalent to the composite of
$\{0\}\times f_0:(0, s_0)\to (0,t_0)$ and $\Delta^1=\Delta^1\times \{t_0\}\hookrightarrow  \Delta^1\times U_A^{op}$.
The functor $\sigma$ carries both morphisms to coCartesian morphisms
so that $f$ maps to a coCartesian morphism.
Finally, we consider the case when $f$ factors through $\{1\}\times \overline{U}_A^{op}$.
As in the previous case, we write $\{1\}\times f_0:(1,s_0)\to (1,t_0)$ for $f$.
Consider the morphism $g_{s_0}$ given by $\Delta^1\times\{s_0\}\hookrightarrow \Delta^1\times \overline{U}_A^{op}$.
By the previous case,
the composite $f\circ g_{s_0}$ and $g_{s_0}$ map to coCartesian morphisms.
It follows that $f$ maps to a coCartesian morphism.
\QED

\begin{Construction}
\label{Hconstruction}
The limits of the restriction
\[
\{i\}\times U_A^{op}\to  \CAlg_A^+\stackrel{\QC}{\longrightarrow} \wCat
\]
induces
$\lim_{\Spec C\in U_A}\QC(C\otimes_AS^1)\to \lim_{\Spec C\in U_A}\QC(C)$.
Using the extension to $\Delta^1\times \overline{U}_A^{op}$ and taking $S^1$-invariants,
we obtain the commutative diagram
\[
\label{diagramFF}
\vcenter{
\xymatrix{
\QC((A\otimes_kA)\otimes_AS^1)^{S^1} \ar[r] \ar[d] & \lim_{\Spec C\in U_A}\QC(C\otimes_AS^1)^{S^1} \ar[d] \\
\QC(A\otimes_kA)^{S^1} \ar[r] &  \lim_{\Spec C\in U_A}\QC(C)^{S^1}
}\tag{3.2}}
\]
Since the forgetful functor
$\Mod(\Mod_A^{S^1})_{\mathcal{H}}\to  \Mod(\Mod_A^{S^1})$
preserves cocartesian morphisms,
it follows from Proposition~\ref{cocarsection} and \cite[3.3.3.2]{HTT} that
$\{0\}\times U_A^{op}\to \Mod(\Mod_A^{S^1})_{\mathcal{H}}\to  \Mod(\Mod_A^{S^1})$ induced by $\sigma$, which is a section
of $\{0\}\times U_A^{op}\to \CAlg_A$, determines an object of $\lim_{\Spec C\in U_A}\QC(C\otimes_AS^1)^{S^1}$. We denote the object by $\tilde{\mathcal{H}}$.
\end{Construction}

\begin{Corollary}
\label{semiorigin}
We regard $\HH_\bullet(\textup{pr}_2^*(\CCC)/A)$
as an object of  $\QC((A\otimes_kA)\otimes_AS^1)^{S^1}$ (cf. Construction~\ref{Hconstruction}).
Then the image of $\HH_\bullet(\textup{pr}_2^*(\CCC)/A)$
in $\lim_{\Spec C\in U_A}\QC(C\otimes_AS^1)^{S^1}$
is naturally equivalent to $\tilde{\mathcal{H}}$.
Here $\Perf_A^\otimes$-module structure on
$\textup{pr}_2^*(\CCC)$ in $\HH_\bullet(\textup{pr}_2^*(\CCC)/A)$
defined to be the restriction of the $\Perf_{A\otimes_kA}^\otimes$-module
$\textup{pr}_2^*(\CCC)$ along
$\Perf_A^\otimes\to \Perf_{A\otimes_kA}^\otimes$ given by $A\simeq A\otimes_kk\to A\otimes_kA$.
\end{Corollary}

\Proof
Apply Proposition~\ref{cocarsection}.
\QED

We consider
\[
\xymatrix{
A\otimes_kS^1 \ar[r] & (A\otimes_kA)\otimes_kS^1 \ar[r]& (A\otimes_kA)\otimes_AS^1
}
\]
in $\CAlg_k^{S^1}$.
By $\otimes_kS^1$ and $\otimes_AS^1$ we mean the tensor with $S^1$ in $\CAlg_k$ and
$\CAlg_A$, respectively.
The first arrow is induced by $k\otimes_kA\to A\otimes_kA$
determined by $k\to A$.
The second arrow is the canonical morphism.
The $A$-module structure of $A\otimes_kA$ in $(A\otimes_kA)\otimes_AS^1$
is given by $A\simeq A\otimes_kk\to A\otimes_kA$ determined  by $k\to A$.
We also note that there exists a canonical equivalence $(A\otimes_kA)\otimes_AS^1\simeq A\otimes_k(A\otimes_kS^1)$
in $\CAlg_A^{S^1}$.

\begin{Proposition}
\label{hochschildbasechange}
\begin{enumerate}
\item 
There exists a canonical equivalence 
$\HH_\bullet(\CCC/k)\otimes_{A\otimes_kS^1}((A\otimes_kA)\otimes_AS^1)\simeq \HH_\bullet(\textup{pr}_2^*(\CCC)/A)$ in $\Mod_{(A\otimes_kA)\otimes_A{S^1}}(\Mod_A^{S^1})$.

\item
There exists a canonical equivalence 
$\HH_\bullet(\CCC/k)\otimes_{A\otimes_kS^1}A \simeq \HH_\bullet(\CCC/A)$ in $\Mod_A^{S^1}$.

\end{enumerate}
\end{Proposition}

\Proof
We first prove (1).
According to \cite[Lemma 6.3]{IMA} there exists
a canonical equivalence
\[
\HH_\bullet(\CCC/k)\otimes_{A\otimes_kS^1}((A\otimes_kA)\otimes_kS^1)\simeq \HH_\bullet(\textup{pr}_2^*(\CCC)/k)
\]
in $\Mod_{(A\otimes_kA)\otimes_kS^1}(\Mod_A^{S^1})$.
By \cite[Theorem 4.13]{DC1}
there exists a canonical equivalence
\[
\HH_\bullet(\textup{pr}_2^*(\CCC)/k)\otimes_{((A\otimes_kA)\otimes_kS^1)}((A\otimes_kA)\otimes_AS^1)\simeq \HH_\bullet(\textup{pr}_2^*(\CCC)/A)
\]
in $\Mod_{((A\otimes_kA)\otimes_AS^1)}(\Mod_A^{S^1})$.
Combining two equivalences we obtain
$\HH_\bullet(\CCC/k)\otimes_{A\otimes_kS^1}((A\otimes_kA)\otimes_AS^1)\simeq \HH_\bullet(\textup{pr}_2^*(\CCC)/A)$.

The assertion (2) is the special case of \cite[Theorem 4.27]{DC1}.
\QED

\section{Modular interpretaion}

\label{modularsec}

Let $\CCC$ be an $A$-linear stable small $\infty$-category (cf. \cite[Section 3]{DC1}).

\subsection{}
We begin by introducing the purpose of Section~\ref{modularsec}.
We let $\mathbb{T}_{A/k}[-1]$ denote the dg Lie algebra that corresponds to $\widehat{S\times_kS}$
via the categorical equivalence $Lie_A\simeq \FST_A$. The underlying complex is equivalent
to the $(-1)$-shifted tangent module $\mathbb{T}_{A/k}$ of $\Spec A$ over $k$.
Let $\mathbb{T}_{A/k}[-1]^{S^1}$ be the dg Lie algebra obtained by cotensor by $S^1\in\SSS$.
The Hochschild homology $\HH_\bullet(\CCC/A)$ admits the action of the dg Lie algebra $\mathbb{T}_{A/k}[-1]^{S^1}$, which we call the canonical $\mathbb{T}_{A/k}[-1]^{S^1}$-module
(cf. \cite[Definition 6.3]{DC1}).
This action is constructed by using the algebraic structure of the Hochschild pair
$(\HH^\bullet(\CCC/A),\HH_\bullet(\CCC/A))$
and the Kodaira-Spencer morphism for $\CCC$.
See \cite{I} or \cite{IMA} for the convention on the Hochschild cohomology $A$-module spectrum $\HH^\bullet(\CCC/A)\in \Mod_A$ and 
the Hochschild homology $A$-module spectrum $\HH_\bullet(\CCC/A))\in \Mod_A$.
The algebraic structure of the Hochshcild pair is described as 
an algebra over the Kontsevich-Soibelman (topological) operad $\KS$.
We shall refer the reader to \cite{I} for the construction.
Let $\GG_{\CCC}$ denote the dg Lie algebra associated to the $\etwo$-algebra $\HH^\bullet(\CCC/A)$.
Using the algebra $(\HH^\bullet(\CCC/A),\HH_\bullet(\CCC/A))$ over $\KS$
we can define a morphism
$\widehat{A}^L_{\CCC}:\GG_{\CCC}^{S^1}\to \End^L(\HH_\bullet(\CCC/A))$ in $\Fun(BS^1,Lie_A)=Lie_A^{S^1}$,
where $\GG_{\CCC}^{S^1}$ is obtained from $\GG_{\CCC}$ by cotensor by $S^1$,
and $\End^L(\HH_\bullet(\CCC/A))$ is the dg Lie algebra endowed with $S^1$-action
which is associated to the endomorphism
algebra obejct $\End(\HH_\bullet(\CCC/A))\in \Alg_1(\Mod_A^{S^1})$ of $\HH_\bullet(\CCC/A)\in \Mod_A^{S^1}$ (see \cite[Section 6.1]{DC1} for the detail).
The Kodaira-Spencer morphism for $\CCC$ is defined to be a morphism 
$KS_{\CCC}:\mathbb{T}_{A/k}[-1]\to \GG_{\CCC}$ of dg Lie algebras (see \cite[Section 5]{DC1}).
If we think of $\mathbb{T}_{A/k}[-1]$ and $\GG_{\CCC}$ as dg Lie algebra equiped with trivial 
$S^1$-actions, 
$KS_{\CCC}$ is naturally promoted to a morphism in $Lie_A^{S^1}$.
The canonical $\mathbb{T}_{A/k}[-1]^{S^1}$-module $\HH_\bullet(\CCC/A)$
is defined by the composite $\mathbb{T}_{A/k}[-1]^{S^1} \stackrel{KS_{\CCC}^{S^1}}{\to} \GG_{\CCC}^{S^1}\to  \End^L(\HH_\bullet(\CCC/A))$ in $Lie_A^{S^1}$.
We will describe the canonical $\mathbb{T}_{A/k}[-1]^{S^1}$-module $\HH_\bullet(\CCC/A)$ in terms of cyclic deformations of $\HH_\bullet(\CCC/A)$ (see Section~\ref{cyclicdeformationsec}, \cite{IMA}).
The presentation will be given in Proposition~\ref{modularmodule}.

\subsection{}
\label{liemodule}
According to \cite[Lemma 7.17]{DC1}, there is a canonical equivalence
\begin{eqnarray*}
R_{\CCC}:\Rep(\mathbb{T}_{A/k}[-1]^{S^1})(\Mod_A^{S^1})\times_{\Mod_A^{S^1}}\{\HH_\bullet(\CCC/A)\} \ \ \ \ \ \ \ \ \ \ \ \ \ \ \ \ \ \ \ \ \ \ \ \ \ \ \  \ \ \ \ \ \ \ \ \ \ \ \ \ \ \ \ \ \ \ \ \ \  \ \\
\ \ \ \ \ \ \ \ \ \ \ \ \ \ \ \ \ \ \ \ \ \ \ \ \ \ \ \ \ \ \ \simeq \lim_{\Spec C\in (\TSZ_A)_{/\widehat{S\times_kS}}} \Rep(\DD_\infty(C\otimes_AS^1))(\Mod_A^{S^1})\times_{\Mod_A^{S^1}}\{\HH_\bullet(\CCC/A)\}
\end{eqnarray*}
(see Section~\ref{formalstacksec}, \cite{H} for $\DD_\infty$).
If we consider the image of the canonical $\mathbb{T}_{A/k}[-1]^{S^1}$-module $\HH_\bullet(\CCC/A)$
in the $\infty$-category on the right side, there is its presentation in terms of cyclic deformations.
Let
$\mathcal{F}^\circlearrowleft_{A\oplus \End(\HH_\bullet(\CCC/A))}:\EXT_A\to \SSS$
be the functor informally defined by 
\[
C\mapsto \bigl(\Rep(\DD_\infty( C\otimes_AS^1))(\Mod_A^{S^1})\times_{\Mod_A^{S^1}}\{\HH_\bullet(\CCC/A)\}\bigr)^{\simeq},
\]
where $\Rep(\DD_\infty(C\otimes_AS^1))(\Mod_A^{S^1})\to \Mod_A^{S^1}$ is the forgetful functor
(see \cite[Section 8.2 and Remark 8.14]{IMA} and \cite[Remark 6.6]{DC1} for the detail).

Let 
$\LMod\circ \DD_1:\Alg_1^+(\Mod_A)\stackrel{\DD_1}{\to} \Alg_1(\Mod_A)^{op} \stackrel{\LMod}{\to} \wCat$
denote the composite
where the first fuctor is given by the $\eone$-Koszul duality functor $\DD_1:\Alg_1^+(\Mod_A)\to \Alg_1^+(\Mod_A)^{op}$ (see e.g. \cite{IMA}, \cite{DC1}) followed by the forgetful functor $\Alg_1^+(\Mod_A)\to \Alg_1(\Mod_A)$ (here we slightly abuse notation),
and the second functor $\LMod$ indicates the functor corresponding to
the Cartesian fibration
$\LMod(\Mod_A)\to \Alg_1(\Mod_A)$
(see \cite[Section 2]{IMA} or \cite[4.2.1]{HA} for the notation).
Note that 
$\Rep:(Lie_A)^{op}\to \wCat$ is the composite of $\LMod:\Alg_1(\Mod_A)^{op}\to \wCat$ and the universal enveloping algebra functor $U_1:Lie_A\to \Alg_1(\Mod_A)$,
and  there exists $U_1\circ \DD_\infty\simeq \DD_1$ between functors
$\EXT_A\to \Alg_1^+(\Mod_A)$ (see \cite[Proposition 3.3]{IMA}).
Thus, there exists a canonical equivalence
$\LMod\circ \DD_1 |_{\EXT_A}\simeq \Rep\circ \DD_\infty|_{\EXT_A}:\EXT_A\to \wCat$
between functors $\EXT_A\to \wCat$.

We briefly review 
 the Koszul duality functor
\[
\QC(C)=\Mod_{C}(\Mod_A) \to \LMod_{\DD_1(C)}(\Mod_A)
\]
for $C \in \CAlg_A^+$, which sends $P$ to $P\otimes_{C}A$.
Here we abuse notation by writing $\QC$ for $\CAlg_A^+\stackrel{\textup{forget}}{\longrightarrow} \CAlg_A\stackrel{\QC}{\to} \wCat$. 
Note that
$\DD_1(C)\otimes_A C \to A$
exhibits $A$ as a $\DD_1(C)$-$C$-bimodule.
This integral kernel $A$ determines a functor
$I_C:\Mod_{C}(\Mod_A) \to \LMod_{\DD_1(C)}(\Mod_A)$
given by $P\mapsto P\otimes_{C}A$.
 If $C\in \EXT_A$, then
$I_C$ is fully faithful (see e.g. \cite[2.3.6]{H}).
By the construction in \cite[Remark 5.9]{IMA},
$I_C$ is functorial in $C\in \CAlg^+_A$. That is,
there is a natural transformation
\[
\mathcal{I}:\QC \to \LMod \circ \DD_1
\]
between functors $\CAlg^+_A \to \wCat$, such that the evaluation at each $R\in \CAlg^+_A$ is equivalent to $I_R$.
This natural transformation (its $S^1$-equivariant version) determines $J_{\HH_\bullet(\CCC/A)}^\circlearrowleft:\Def^\circlearrowleft(\HH_\bullet(\CCC/A))\to \mathcal{F}^\circlearrowleft_{A\oplus \End(\HH_\bullet(\CCC/A))}$ 
such that for $C\in \EXT_A$, the evaluation at $C$ is 
\begin{eqnarray*}
\bigl(\Mod_{C\otimes_AS^1}(\Mod_A^{S^1})\times_{\Mod_A^{S^1}}\{\HH_\bullet(\CCC/A)\}\bigr)^{\simeq}&\to& \bigl(\LMod_{\DD_1(C\otimes_AS^1)}(\Mod_A^{S^1})\times_{\Mod_A^{S^1}}\{\HH_\bullet(\CCC/A)\}\bigr)^{\simeq}  \\
 &\simeq& \bigl(\Rep(\DD_\infty( C\otimes_AS^1))(\Mod_A^{S^1})\times_{\Mod_A^{S^1}}\{\HH_\bullet(\CCC/A)\}\bigr)^{\simeq}
\end{eqnarray*}
induced by $I_{C\otimes_AS^1}$
(see \cite[Proposition 7.1]{IMA} for detail).

\subsection{}
\label{DE}
Let us formulate Proposition~\ref{modularmodule}.
Let
$\sigma_0:\{0\}\times U_A^{op} \hookrightarrow \Delta^1 \times\overline{U}_A^{op}\to\Mod(\Mod_A^{S^1})_{\mathcal{H}}$
be the restriction of the functor $\sigma$ in Proposition~\ref{cocarsection}.
Since $\sigma_{0}$ carries any morphism to a coCartesian morphism,
it gives rise to an object of
\[
\lim_{ f:\Spec C\to S\times_kS\in U_A}\Def^\circlearrowleft(\HH_\bullet(\CCC/A))(C)
\]
(cf. \cite[3.3.3.2]{HTT}).
We shall write $D_{\CCC}^\circlearrowleft$ for it.
Let $\textup{pr}_2^*(\CCC)$ be the base change $\CCC\otimes_{\Perf_A}\Perf_{A\otimes_kA}\in \ST_{A\otimes_kA}$
of $\CCC\in\ST_A$ (cf. Proposition~\ref{cocarsection}).
For $f:\Spec C\to S\times_kS$ in $U_A$,
we write $\CCC_f$ for $f^*(\textup{pr}_2^*(\CCC))=\textup{pr}_2^*(\CCC)\otimes_{A\otimes_kA}C$.
Then $D_{\CCC}^\circlearrowleft$ is informally described as
the homotopy coherent diagram of the collection
\[
\{(\HH_\bullet(\CCC_f/A), \HH_\bullet(\CCC_f/A)\otimes_{C\otimes_AS^1}A\simeq \HH_\bullet(\CCC/A))\}_{f\in  U_A}.
\]
By definition, the image of $D_{\CCC}^{\circlearrowleft}$ in 
$\lim_{\Spec C\to S\times S\in U_A}\QC(C\otimes_AS^1)^{S^1}$
is $\tilde{\mathcal{H}}$ (cf. Construction~\ref{Hconstruction}).

Passing to limits, the natural transformation
$J_{\HH_\bullet(\CCC/A)}^\circlearrowleft:\Def^\circlearrowleft(\HH_\bullet(\CCC/A))\to \mathcal{F}^\circlearrowleft_{A\oplus \End(\HH_\bullet(\CCC/A))}$ 
sends $D_{\CCC}^\circlearrowleft$ to an object of
$\lim_{ f:\Spec C\to S\times_kS\in U_A}\mathcal{F}^\circlearrowleft_{A\oplus \End(\HH_\bullet(\CCC/A))}(C)$ which we denote by $E_{\CCC}^{\circlearrowleft}$.
We remark that ``$\lim_{f:\Spec C\to S\times_kS\in U_A}$'' can be replaced
with ``$\lim_{f:\Spec C\to \widehat{S\times_kS}\in (\TSZ_A)_{/\widehat{S\times_kS}}}$''
(cf. Lemma~\ref{nonformalcompletion}).

\begin{Proposition}
\label{modularmodule}
The canonical $\mathbb{T}_{A/k}[-1]^{S^1}$-module $\HH_\bullet(\CCC/A)$
corresponds to $E_{\CCC}^\circlearrowleft$ through the equivalence $R_{\CCC}$.
\end{Proposition}

\begin{Remark}
Note the underlying object of $D_{\CCC}^\circlearrowleft$ is
$\tilde{\mathcal{H}}$, and $D_{\CCC}^\circlearrowleft$ maps to $E_{\CCC}^\circlearrowleft$.
Roughly, Proposition~\ref{modularmodule} means that
there is a recipe to obtain 
the canonical $\mathbb{T}_{A/k}[-1]^{S^1}$-module $\HH_\bullet(\CCC/A)$
from cyclic deformations of $\HH_\bullet(\CCC/A)$ arising from deformations of $\CCC$.
\end{Remark}

\Proof
According to \cite[Lemma 8.13]{IMA},
for $M\in Lie_A$ there exists a canonical equivalence
\begin{eqnarray*}
\Map_{\Fun(\EXT_A,\SSS)}(\mathcal{F}_{M},\mathcal{F}^\circlearrowleft_{A\oplus \End(\HH_\bullet(\CCC/A))}) &\simeq&  \Rep(M^{S^1})(\Mod_A^{S^1})^\simeq\times_{(\Mod_A^{S^1})^{\simeq}}\{\HH_\bullet(\CCC/A)\} \\
\end{eqnarray*}
which is functorial in $M\in Lie_A$. Here $\mathcal{F}_M$ is a pointed formal stack associated to $M$.
Through this equivalence for $M=\mathbb{T}_{A/k}[-1]$
and $R_{\CCC}$, $E_{\CCC}^{\circlearrowleft}$  is classified,
as the object (on the right side), by the following composite of $\overline{KS}_{\CCC}$ and maps
in Theorem 1.2 in \cite{IMA}:
\[
u:\mathcal{F}_{\mathbb{T}_{A/k}[-1]}\simeq \widehat{S\times_kS} \stackrel{\overline{KS}_\CCC}{\longrightarrow}  \Def_{\CCC}\stackrel{M_\CCC^{\circlearrowleft}}{\longrightarrow} \Def^{\circlearrowleft}(\HH_\bullet(\CCC/A))\stackrel{J_{\HH_\bullet(\CCC/A)}^{\circlearrowleft}}{\longrightarrow} \mathcal{F}^\circlearrowleft_{A\oplus \End(\HH_\bullet(\CCC/A))}
\]
(this is an obvious consequence of the definition of the above sequence in \cite{IMA}).
For $X\in \FST_A$, there exists a canonical equivalence
\begin{eqnarray*}
\Map_{\Fun(\EXT_A,\SSS)}(X,\Def^{\circlearrowleft}(\HH_\bullet(\CCC/A))) &\simeq&  (\QC_H^{\circlearrowleft}(X)^{S^1})^\simeq\times_{(\Mod_A^{S^1})^{\simeq}}\{\HH_\bullet(\CCC/A)\} \\
\end{eqnarray*}
which is functorial in $X\in \FST_A$ (see Section~\ref{cyclicdeformationsec}).
In a similar vein, $D_{\CCC}^\circlearrowleft$ is the object of
\[
\Map_{\Fun(\EXT_A,\SSS)}(\widehat{S\times_kS},\Def^{\circlearrowleft}(\HH_\bullet(\CCC/A)))\simeq \lim_{ \Spec C\in (\TSZ_A)_{/\widehat{S\times_kS}}}\Def^\circlearrowleft(\HH_\bullet(\CCC/A))(C)
\]
that corresponds to 
$\widehat{S\times_kS} \stackrel{\overline{KS}_\CCC}{\longrightarrow}  \Def_{\CCC}\stackrel{M_\CCC^{\circlearrowleft}}{\longrightarrow} \Def^{\circlearrowleft}(\HH_\bullet(\CCC/A))$
through the equivalence (this is an obvious consequence of the defintions
of $\overline{KS}_\CCC$ and $M_\CCC^{\circlearrowleft}$ in \cite{IMA}).

Now our claim follows from \cite[Theorem 1.2]{IMA}, which says that $u$
is equivalent to
\[
v:\mathcal{F}_{\mathbb{T}_{A/k}[-1]}\simeq \widehat{S\times_kS} \stackrel{\overline{KS}_\CCC}{\longrightarrow}  \Def_{\CCC}\stackrel{J_\CCC}{\longrightarrow} \FF_{\GG_{\CCC}}\stackrel{}{\longrightarrow} \mathcal{F}^\circlearrowleft_{A\oplus \End(\HH_\bullet(\CCC/A))}.
\]
This sequence appears in \cite[Theorem 1.2]{IMA}, and
the final arrow is defined in \cite[Construction 8.6]{IMA}
(in {\it loc. cit.}, we denote it by $\mathfrak{T}^{\etwo}_{\CCC}:\mathcal{F}_{A\oplus \HH^\bullet(\CCC/A)}\to \mathfrak{h}_{A\oplus \End(\HH_\bullet(\CCC/A))}\circ h\circ \DD_2 \simeq \mathfrak{h}_{A\oplus \End(\HH_\bullet(\CCC/A))}\circ \DD_1\circ h$). The final arrow corresponds to an object of 
\[
\Rep(\GG_{\CCC}^{S^1})(\Mod_A^{S^1})\times_{\Mod_A^{S^1}}\{\HH_\bullet(\CCC/A)\}
\]
determined by
$\widehat{A}^L_{\CCC}:\GG_{\CCC}^{S^1} \to \End^L(\HH_\bullet(\CCC/A))$.
By definition (\cite[Section 5]{DC1}), the Kodaira-Spencer morphism $\mathbb{T}_{A/k}[-1]\to \GG_{\CCC}$
corresponds to the composite
$\mathcal{F}_{\mathbb{T}_{A/k}[-1]}\simeq \widehat{S\times_kS} \to  \Def_{\CCC}\to \FF_{\GG_{\CCC}}$
via $Lie_A\simeq \FST_A$.
We deduce that $v$ is classified by 
an object of 
$\Rep(\mathbb{T}_{A/k}[-1]^{S^1})(\Mod_A^{S^1}) \times_{\Mod_A^{S^1}}\{\HH_\bullet(\CCC/A)\}$ determined by 
the composite
$\mathbb{T}_{A/k}[-1]^{S^1}\stackrel{KS_{\CCC}^{S^1}}{\to} \GG_{\CCC}^{S^1}\stackrel{\widehat{A}^L_{\CCC}}{\to}  \End^L(\HH_\bullet(\CCC/A))$, that is, the canonical action of $\mathbb{T}_{A/k}[-1]^{S^1}$
on $\HH_\bullet(\CCC/A)$.
\QED

\section{Quasicoherent complexes between Ind-coherent complexes
and Lie algebra modules}

In this section, for the reader's convenience, we review the diagram \cite[Section 7.3]{DC1}, which  involves the $\infty$-cateory of modules over a dg Lie algebra
and the $\infty$-category of Ind-coherent complexes over a formal stack,
see 
Proposition~\ref{diagramprop}.

\label{importantpreparesec}

\subsection{}
\label{zigzagconstructionsec}

By abuse of notation we continue to write 
$\QC:\CAlg_A^+\to \wCat$
for the composite functor $\CAlg_A^+\stackrel{\textup{forget}}\to \CAlg_A\to \wCat$
where the second functor corresponds to the coCartesian fibration
$\Mod(\Mod_A)\to \CAlg_A$.
Namely, $\QC$ carries $C\in \CAlg_A^+$ to $\Mod_C$,
and we write $\QC(C)$ for $\Mod_C$.
Let $\QC_!|_{\EXT_A}:\EXT_A\to \wCat$
be the functor given on objects by $C\mapsto \QC_!(\Spec C)=\Ind(\Coh(C))$.
Here $\Coh(C)$ is the full subcategory of $\QC(C)$ spanned by
those objects which are 
bounded with coherent cohomology (with respect to the standard $t$-structure).
For $f:\Spec C'\to \Spec C$ in $\TSZ_A$, it carries $f$ to
the $!$-pullback functor $f^!:\QC_!(\Spec C)\to \QC_!(\Spec C')$
which is the right adjoint to the proper pushforward functor 
 $f^{\textup{IndCoh}}_*=\Ind(f_*|_{\Coh(\Spec C')}):\QC_!(\Spec C')=\Ind(\Coh(\Spec C'))\to \QC_!(\Spec C)=\Ind(\Coh(\Spec C))$.
This functor is the restriction of the functor $\QC_!:\FST_!\to \wCat$
constructed in \cite{Ind}, \cite{Gai2} (see also \cite{DC1}).

We consider three functors from $\EXT_A$ to $\wCat$: $\QC|_{\EXT_A},\ \QC_!|_{\EXT_A},\ \LMod\circ \DD_1 |_{\EXT_A}$.
Taking right Kan extensions of $\QC|_{\EXT_A},\ \QC_!|_{\EXT_A},\ \LMod\circ \DD_1 |_{\EXT_A}$ along $\EXT_A=(\TSZ_A)^{op}\to (\GFST_A)^{op}$,
we define three functors
\[
\QC_H',\ \QC_!',\ \Rep_H' :(\GFST_A)^{op}\longrightarrow \wCat.
\]

We let $\Upsilon|_{\EXT_A}:\QC|_{\EXT_A}\to \QC_!|_{\EXT_A}$ denote the natural transformation
between functors $\EXT_A\to \wCat$, which is
induced by $\Upsilon$ (cf. Section~\ref{NC}).
Recall $\mathcal{I}|_{\EXT_A}:\QC|_{\EXT_A} \to \LMod \circ \DD_1|_{\EXT_A}$ from Section~\ref{liemodule}.
Let
\[
\QC_!' \stackrel{\Upsilon'}{\longleftarrow} \QC_H' \stackrel{\mathcal{I}'}{\longrightarrow} \Rep_H'
\]
be the diagram obtained from 
$\QC_!|_{\EXT_A} \stackrel{\Upsilon|_{\EXT_A}}{\longleftarrow}\QC|_{\EXT_A} \stackrel{\mathcal{I}|_{\EXT_A}}{\longrightarrow} \LMod \circ \DD_1|_{\EXT_A}$
by taking the right Kan extensions.
For $W\in \GFST_A$, 
$\QC_!'(W) \leftarrow \QC_H'(W) \to \Rep_H'(W)$
is naturally equivalent to
\[
\lim_{\Spec C\in (\TSZ_A)_{/W}}\QC_!(\Spec C)\leftarrow \lim_{\Spec C\in (\TSZ_A)_{/W}}\QC(\Spec C)\to\lim_{\Spec C\in (\TSZ_A)_{/W}}\LMod_{\DD_1(C)}.
\]
Both $\QC_H'(W)\rightarrow \QC_!'(W)$
and $\QC_H'(W) \to \Rep_H'(W)$ are fully faithful.

For $\Spec B\in (\Aff_k)_{S//S}$
we let $\widehat{\Spec B}$ denote the functor $\EXT_A\to \SSS$
defined as the restriction of the functor $(\CCAlg_k)_{A//A}\to \SSS$ corepresented by $B$
(cf. \cite[Definition 2.2.7]{H}, \cite[Section 3]{DC1}, Section~\ref{formalstacksec}).
According to \cite[2.2.8]{H}, $\widehat{\Spec B}$ lies in $\FST_A\subset \Fun(\EXT_A,\SSS)$.
Let $\textup{comp}:(\Aff_k)_{S//S}\to \FST_A$ be the formal completion functor
given by the assignment $\Spec B\mapsto \widehat{\Spec B}$.
The composite $(\Aff_k)_{S//S}\to \FST_A\simeq Lie_A$ is naturally equivalent to
the functor $(\Aff_k)_{S//S}\subset (\CAlg_A^+)^{op}\stackrel{\DD_\infty}{\longrightarrow} Lie_A$.
By using the definition of the formal completion, we easily see:

\begin{Lemma}
\label{nonformalcompletion}
Let $Y$ be an object of $(\Aff_k)_{S//S}$.
The functor $\textup{comp}$ induces an equivalence of $\infty$-categories
$\bigl((\Aff_k)_{S//S}\bigr)_{/Y}\times_{ (\Aff_k)_{S//S}}\TSZ_A\stackrel{\sim}{\to} (\TSZ_A)_{/\widehat{Y}}$.
\end{Lemma}

Given $M\in \SSS$,  we let $\otimes_AM:(\CCAlg_A)^+:=(\CCAlg_k)_{A//A} \to (\CCAlg_A)^+$ 
denote the functor given by tensor by $M$ in $(\CCAlg_A)^+$.
Consider the composite
\[
\xi_M: (\CCAlg_A)^+ \stackrel{\otimes_AM}{\longrightarrow} (\CCAlg_A)^+\stackrel{\textup{comp}^{op}}{\longrightarrow} (\FST_A)^{op}\stackrel{\Theta_A^{op}}{\longrightarrow}(\GFST_A)^{op}.
\]
We write $(\Spec C)^\wedge$ for the image of $\widehat{\Spec C}$
under the equivalence $\Theta_A:\FST_A\stackrel{\sim}{\to} \GFST_A$.
The composite $\xi_M$ carries $B$ to $(\Spec B\otimes_AM)^\wedge$.

Let $X$ be a pointed formal stack over $A$, that is, an object of $\FST_A$.
We set $(\TSZ_A)_{/X}=\TSZ_A\times_{\FST_A}(\FST_A)_{/X} $
and consider the composite
\[
\eta_{M,G}:\bigl((\TSZ_A)_{/X}\bigr)^{op}\stackrel{\textup{forget}}{\longrightarrow} (\CCAlg_A)^+ \stackrel{\xi_M}{\longrightarrow} (\GFST_A)^{op}\stackrel{G}{\to} \wCat
\]
where the third functor $G$ is either 
$\QC_H',\ \QC_!'$ or  $\Rep_H'$. Define
\[
\xymatrix{
\QC_!^{M\wedge}(X) & 
\QC_H^{M\wedge}(X)  \ar[r] \ar[l]  & \Rep_H^{M\wedge}(X)  \\
}
\] 
to be the diagram obtained from $\QC_!' \leftarrow \QC_H' \to \Rep_H'$ by taking the limits of $\eta_{M,G}$ in $\wCat$.
We think of the diagram as an object of $\Fun(\Delta^1\sqcup_{\{0\}}\Delta^1,\wCat)$
and denote it by $\mathbf{D}_{M,X}$. 
The diagram $\mathbf{D}_{M,X}$ can naturally be identified with
\[
\xymatrix{
\lim_{\Spec C\in (\TSZ_A)_{/X}}\QC_H'((\Spec C\otimes_AM)^\wedge) \ar[r] \ar[d] & \lim_{\Spec C\in (\TSZ_A)_{/X}}\Rep_H'((\Spec C\otimes_AM)^\wedge) \\
\lim_{\Spec C\in (\TSZ_A)_{/X}}\QC'_!((\Spec C\otimes_AM)^\wedge).
}
\]
These two functors are obviously fully faithful.

Next we consider the cases when $M=S^1$ and $M=\ast$.
We set
$\QC_!^{\circlearrowleft\wedge}(X)=\QC_!^{S^1\wedge}(X),\ \QC_H^{\circlearrowleft\wedge}(X)=\QC_H^{S^1\wedge}(X),\ \Rep_H^{\circlearrowleft\wedge}(X)=\Rep_H^{S^1\wedge}(X),\ \QC_!^{\wedge}(X)=\QC_!^{\ast\wedge}(X),\ \QC_H^{\wedge}(X)=\QC_H^{\ast\wedge}(X),\ \Rep_H^{\wedge}(X)=\Rep_H^{\ast\wedge}(X)$.
The construction of $\mathbf{D}_{M,X}$ is functorial with respect to $M$:
the assignment $M \mapsto \mathbf{D}_{M,X}$ can be promoted to $\SSS\to \Fun(\Delta^1\sqcup_{\{0\}}\Delta^1,\wCat)$.
In particular, the $S^1$-equivariant map $S^1\to \ast$
induces a morphism $\mathbf{D}_{S^1,X}\to \mathbf{D}_{\ast,X}$ in $\Fun(BS^1\times  (\Delta^1\sqcup_{\{0\}}\Delta^1),\wCat)$ (see \cite[Remark 7.11, Construction 7.12]{DC1} for the formulation).
Furthermore, we focus on the case when $X=\widehat{S\times_kS}$.
Combined with these observations,
the following is proved in \cite[Proposition 7.13, Proposition 7.14, Construction 7.15]{DC1}:

\begin{Proposition}[\cite{DC1}]
\label{diagramprop}
The followings hold:
\begin{enumerate}

\item 
There exists the diagram
\[
\xymatrix{
\ \ \ \ \ \QC_!((S\times_kLS)_S^\wedge)  \simeq \QC_!^{\circlearrowleft\wedge}(\widehat{S\times_kS})
 \ar[d] & \QC_H^{\circlearrowleft\wedge}(\widehat{S\times_kS}) \ar[r]^(0.3){\beta'} \ar[l]_(0.3)\beta \ar[d] & \Rep_H^{\circlearrowleft\wedge}(\widehat{S\times_kS})\simeq \Rep(\mathbb{T}_{A/k}[-1]^{S^1})(\Mod_A) \ar[d] \\ 
\ \ \ \ \  \QC_!((S\times_kS)_S^\wedge)   \simeq  \QC_!^{\wedge}(\widehat{S\times_kS})
 & \QC_H^{\wedge}(\widehat{S\times_kS}) \ar[r] \ar[l] & \Rep_H^{\wedge}(\widehat{S\times_kS})\simeq \Rep(\mathbb{T}_{A/k}[-1])(\Mod_A) 
}
\]
in $\Fun(BS^1,\wCat)$. This diagram up to equivalences is $\mathbf{D}_{S^1,X}\to \mathbf{D}_{\ast,X}$.
The vertical functor on the right side is
determined by the restriction along the diagonal morphism
$\mathbb{T}_{A/k}[-1]\to\mathbb{T}_{A/k}[-1]^{S^1}$.
The vertical functor on the left side is the $!$-pullback functor
along $(\textup{id}_S\times\iota)_S^\wedge:(S\times_kS)_S^\wedge\to (S\times_kLS)_S^\wedge$.
Every horizontal functors is fully faithful.

\item 
Taking $S^1$-invariants, we obtain
\[
\xymatrix{
\QC_!((S\times_kLS)_S^\wedge)^{S^1}
 \ar[d] & \QC_H^{\circlearrowleft\wedge}(\widehat{S\times_kS})^{S^1} \ar[r]^{(\beta')^{S^1}} \ar[l]_(0.4){\beta^{S^1}} \ar[d] & \Rep(\mathbb{T}_{A/k}[-1]^{S^1})(\Mod_A^{S^1}) \ar[d] \\ 
 \QC_!((S\times_kS)_S^\wedge)^{S^1}
 & \QC_H^{\wedge}(\widehat{S\times_kS})^{S^1} \ar[r] \ar[l] & \Rep(\mathbb{T}_{A/k}[-1])(\Mod_A^{S^1}) 
}
\]
where horizontal functors are fully faithful functors.
\end{enumerate}
\end{Proposition}

\section{Formal completion}

\label{completion}

Let $W$ be a pointed formal stack over $A$.
By definition (see Defintion~\ref{QCHcircle}), 
$\QC_H^{\circlearrowleft}(W)$ is equivalent to $\lim_{\Spec C\in (\TSZ_A)_{/W}}\QC(C\otimes_AS^1)$.
If we remember the theory of formal schemes,
an object of $\QC(C\otimes_AS^1)=\Mod_{C\otimes_AS^1}(\Mod_A)$ can not be thought of as
being formally complete along $\Spec A\to \Spec C\otimes_AS^1$.
In this section, we construct a sort of formal completions.
We start with a general situation.

Let us consider the sequence
\[
\EXT_A\hookrightarrow \bigl((\Aff_k)_{S//S}\bigr)^{op}\stackrel{\textup{comp}^{op}}{\longrightarrow} (\FST_A)^{op}\stackrel{\Theta_A^{op}}{\simeq} (\GFST_A)^{op}
\]
where the first functor is the evident inclusion.
Then this sequence induces the adjoint pairs
\[
\xymatrix{
\Fun( (\GFST_A)^{op},\wCat) \ar@<0.5ex>[r]^(0.45){\textup{res}_2}  &  \Fun((\CCAlg_k)_{A//A},\wCat)  \ar@<0.5ex>[r]^(0.55){\textup{res}_1} \ar@<0.5ex>[l]^(0.55){\mathcal{R}_2} &  \Fun(\EXT_A,\wCat) \ar@<0.5ex>[l]^(0.4){\mathcal{R}_1}.
}
\]
The left adjoint functors are given by restrictions.
The right adjoint $\mathcal{R}_1$ and $\mathcal{R}_2$
are given by right Kan extensions along $\EXT_A\to (\CCAlg_k)_{A//A} $ and
$ (\CCAlg_k)_{A//A} \to (\GFST_A)^{op}$, respectively.

\begin{Construction}
\label{abstractcompletion}
Let $F:(\CCAlg_k)_{A//A}\to \wCat$ be a functor. The typical example is
$\QC:(\CCAlg_k)_{A//A}\to \wCat$ which carries $[A\to B\to A]$ to $\Mod_B$.
Let $p:F\to \mathcal{R}_1\circ \textup{res}_1(F)$ be the unit map determined by the adjoint pair $(\textup{res}_1,\mathcal{R}_1)$.
Let $q:\textup{res}_2\circ \mathcal{R}_2(\mathcal{R}_1\circ \textup{res}_1(F))\to \mathcal{R}_1\circ \textup{res}_1(F)$ be the counit map determined by the adjoint pair
$(\textup{res}_2,\mathcal{R}_2)$. We obtain the diagram
\[
F\stackrel{p}{\to}
 \mathcal{R}_1\circ \textup{res}_1(F) \stackrel{q}{\leftarrow} \textup{res}_2\circ \mathcal{R}_2(\mathcal{R}_1\circ \textup{res}_1(F)).
\]
Observe that $q$ is an equivalence.
Note that $\mathcal{R}_2(\mathcal{R}_1\circ \textup{res}_1(F))$ is a right Kan extension of
$F|_{\EXT_A}=\textup{res}_1(F)$ along the fully faithful embedding $\EXT_A\to (\GFST_A)^{op}$.
Thus, for $B\in (\CCAlg_k)_{A//A}$, $\textup{res}_2\circ \mathcal{R}_2(\mathcal{R}_1\circ \textup{res}_1(F))(B)$
is described as the limit $\lim_{\Spec C\in (\TSZ_A)_{/(\Spec B)^\wedge}}F(C)$.
Similarly,  $(\mathcal{R}_1\circ \textup{res}_1(F))(B)$ is $\lim_{\Spec C\in (\TSZ_A)_{/\Spec B}}F(C)$,  where $(\TSZ_A)_{/\Spec B}=\TSZ_A\times_{(\Aff_k)_{S//S}}\bigr((\Aff_k)_{S//S}\bigl)_{/\Spec B}$.
Then $q$ is given by $\lim_{\Spec C\in (\TSZ_A)_{/(\Spec B)^\wedge}}F(C)\to \lim_{\Spec C\in (\TSZ_A)_{/\Spec B}}F(C)$ induced by the equivalence $(\TSZ_A)_{/\Spec B}\simeq (\TSZ_A)_{/(\Spec B)^\wedge}$ which is determined by the 
fully faithful functor $(\TSZ_A)_{/\Spec B}\to (\FST_A)_{/\widehat{\Spec B}}\simeq (\GFST_A)_{/(\Spec B)^\wedge}$, where the first functor is induced by $\textup{comp}$ (the fully faithfulness follows immediately from  the definition of $\widehat{\Spec B}$, see Lemma~\ref{nonformalcompletion}).
It follows that $q$ is an equivalence.
Set $\mathcal{R}=\mathcal{R}_2\circ \mathcal{R}_1$.
Using the inverse $q^{-1}$ of $q$, we obtain 
\[
\textup{comp}_{F}:F\to  \mathcal{R}_1\circ \textup{res}_1(F) \simeq \textup{res}_2\circ \mathcal{R}_2(\mathcal{R}_1\circ \textup{res}_1(F))=\textup{res}_2\circ \mathcal{R}\circ \textup{res}_1(F).
\]

\end{Construction}

\begin{Remark}
In Construction~\ref{abstractcompletion}, for $[A\to B\to A]\in (\CCAlg_k)_{A//A}$,
objects of $\mathcal{R}_1\circ \textup{res}_1(F)(B)\simeq \textup{res}_2\circ \mathcal{R}\circ \textup{res}_1(F)(B)$
should be thought of as ``objects on the formal neighborhood of $\Spec A\to \Spec B$''.
The natural transformation
$\textup{comp}_{F}:F\to  \mathcal{R}_1\circ \textup{res}_1(F) \simeq\textup{res}_2\circ \mathcal{R}\circ \textup{res}_1(F)$ sends an object of $F(B)$ to ``its formal completion along $\Spec A\to \Spec B$''.
\end{Remark}

\begin{Example}
Let $\QC$ denote the composite functor $(\CCAlg_k)_{A//A}\to (\CCAlg_k)_{A/}\hookrightarrow \CAlg_A\to \wCat$ induced by the functor $\QC:\CAlg_A\to \wCat$.
Applying Construction~\ref{abstractcompletion} to $F=\QC$ we define
$\textup{comp}_{\QC}: \QC\to \textup{res}_2\circ \mathcal{R}\circ \textup{res}_1(\QC)=\QC_H'|_{(\CCAlg_k)_{A//A}}$. Here we note that by definition
$\mathcal{R}\circ \textup{res}_1(\QC)=\QC_H'$.
\end{Example}

\begin{Example}
Let $\LMod\circ \DD_1$ denote the composite functor $(\CCAlg_k)_{A//A}\stackrel{\textup{forget}}{\longrightarrow} \Alg_1^+(\Mod_A)\stackrel{\DD_1}{\to} \Alg_1^+(\Mod_A)^{op} \stackrel{\LMod}{\to} \wCat$, that is, the functor which carries
$B$ to $\LMod_{\DD_1(B)}(\Mod_A)$.
Applying Construction~\ref{abstractcompletion} to $F=\LMod\circ \DD_1$ we define
$\textup{comp}_{\LMod\circ \DD_1}: \LMod\circ \DD_1\to \textup{res}_2\circ \mathcal{R}\circ \textup{res}_1(\LMod\circ \DD_1)=\Rep_H'|_{(\CCAlg_k)_{A//A}}$. By definition,
$\mathcal{R}\circ \textup{res}_1(\LMod\circ \DD_1)=\Rep_H'$.
\end{Example}

\begin{Example}
Let $\Rep\circ \DD_\infty$ denote the composite functor $\LMod\circ U_1\circ \DD_\infty:(\CCAlg_k)_{A//A}\to \wCat$, that is, the functor which carries
$B$ to $\Rep(\DD_\infty(B))(\Mod_A)$.
Applying Construction~\ref{abstractcompletion} to $F=\Rep\circ \DD_\infty$ we define
$\textup{comp}_{\Rep\circ \DD_\infty}: \Rep\circ \DD_\infty\to \textup{res}_2\circ \mathcal{R}\circ \textup{res}_1(\Rep\circ \DD_\infty)$. It is possible to prove that $\textup{comp}_{\Rep\circ \DD_\infty}$ is an equivalence
(cf. \cite[Lemma 7.17]{DC1}).
\end{Example}

Let $\mathcal{I}:\QC\to \LMod\circ \DD_1$ be the natural transformation (see Section~\ref{liemodule}).
This functor induces the natural transfromation $\mathcal{I}':\QC_H'\to \Rep_H'$ (see Section~\ref{importantpreparesec}).
Applying Construction~\ref{abstractcompletion} to $\QC\to \LMod\circ \DD_1$, we have
the commutative diagram
\[
\xymatrix{
\QC     \ar[r]^{\mathcal{I}} \ar[d]^{\textup{comp}_{\QC}} &    \LMod\circ \DD_1 \ar[d]^{\textup{comp}_{ \LMod\circ \DD_1 }} \\
\QC_H'|_{(\CCAlg_k)_{A//A}} \ar[r]^{\mathcal{I'}} & \Rep_H'|_{(\CCAlg_k)_{A//A}}
}
\]
in $\Fun((\CCAlg_k)_{A//A},\wCat)$.
If we put $T:\EXT_A\hookrightarrow (\CCAlg_k)_{A//A}   \stackrel{\otimes_AS^1}{\longrightarrow}  \Fun(BS^1,(\CCAlg_k)_{A//A})$ where the second functor is induced by the tensor by $S^1$
in $(\CCAlg_k)_{A//A}$, the composition with $T$
gives rise to the commutative diagram 
\[
\xymatrix{
\QC\circ T     \ar[r] \ar[d]^{\textup{comp}_{\QC}\circ T} &    \LMod\circ \DD_1\circ T \ar[d]^{\textup{comp}_{ \LMod\circ \DD_1}\circ T} \\
\QC_H'|_{(\CCAlg_k)_{A//A}} \circ T\ar[r] & \Rep_H'|_{(\CCAlg_k)_{A//A}}\circ T
}
\]
in $\Fun(\EXT_A,\Fun(BS^1,\wCat))$.

\begin{Proposition}
The natural transformation
$\textup{comp}_{ \LMod\circ \DD_1}\circ T:\LMod\circ \DD_1\circ T \to \Rep_H'|_{(\CCAlg_k)_{A//A}}\circ T$
is an equivalence.
\end{Proposition}

\Proof
This is essentially proved in \cite{DC1}: we here review it.
By definition, for each $C\in \EXT_A$, the induced functor
$\LMod\circ \DD_1\circ T(C) \to \Rep_H'|_{(\CCAlg_k)_{A//A}}\circ T(C)$ (in $\wCat$)
can naturally be identified with 
\begin{eqnarray*}
\Rep(\DD_\infty(C\otimes_AS^1))\simeq \LMod_{\DD_1(C\otimes_AS^1)}(\Mod_A)\to \ \ \ \ \ \ \ \ \ \ \ \ \ \ \ \ \ \ \ \ \ \ \ \ \ \ \ \ \ \ \ \ \ \ \ \ \ \ \ \ \ \ \ \ \ \ \ \ \ \ \ \ \ \ \ \ \ \ \ \ \ \ \  \ \ \ \ \ \ \ \ \ \ \ \\
\ \ \ \ \ \ \ \ \ \ \ \lim_{\Spec R\in (\TSZ_A)_{/(\Spec C\otimes_AS^1)^\wedge}} \LMod_{\DD_1(R)}(\Mod_A)\simeq\lim_{\Spec R\in (\TSZ_A)_{/(\Spec C\otimes_AS^1)^\wedge}} \Rep(\DD_\infty(R))(\Mod_A).
\end{eqnarray*}
where the equivalences in the sequence come from \cite[Proposition 3.3, Proposition 7.1]{IMA}. 
This is an equivalence by \cite[Lemma 7.17]{DC1}. This completes the proof.
\QED

Taking into account the inverse of
$\textup{comp}_{ \LMod\circ \DD_1}\circ T$, we see:

\begin{Corollary}
\label{factorizationcomplete}
The natural transformation
$\QC\circ T     \to    \LMod\circ \DD_1\circ T$
factors as
\[
\QC\circ T \to \QC_H'|_{(\CCAlg_k)_{A//A}} \circ T\to\Rep_H'|_{(\CCAlg_k)_{A//A}}\circ T\simeq \LMod\circ \DD_1\circ T.
\]
\end{Corollary}

\begin{Corollary}
Let $\mathcal{V}_L'\in \Rep(\mathbb{T}_{A/k}[-1]^{S^1})(\Mod_A^{S^1})$
denote the canonical $\mathcal{T}_{A/k}[-1]^{S^1}$-module $\HH_\bullet(\CCC/A)$
(cf. \cite[Definition 6.3, Construction 6.10]{DC1}, Section~\ref{modularsec}).
Then  $\mathcal{V}_L'$ lies in the essential image of the fully faithful functor
$(\beta')^{S^1}:\QC_H^{\circlearrowleft\wedge}(\widehat{S\times_kS})^{S^1} \to \Rep_H^{\circlearrowleft\wedge}(\widehat{S\times_kS})^{S^1}\simeq \Rep(\mathbb{T}_{A/k}[-1]^{S^1})(\Mod_A^{S^1})$.
\end{Corollary}

\Proof
According to Proposition~\ref{modularmodule}, $\mathcal{V}_L'$ is 
the image of $D_{\CCC}^{\circlearrowleft}$ under
\begin{eqnarray*}
 \lim_{\Spec C\in (\TSZ_A)_{/\widehat{S\times S}}}\QC(C\otimes_AS^1)^{S^1}\times_{\Mod_A^{S^1}}\{\HH_\bullet(\CCC/A)\} &\stackrel{\textup{forget}}{\longrightarrow}& \lim_{\Spec C\in (\TSZ_A)_{/\widehat{S\times S}}}\QC(C\otimes_AS^1)^{S^1} \\
&\to&  \lim_{\Spec C\in (\TSZ_A)_{/\widehat{S\times S}}}\LMod_{\DD_1(C\otimes_AS^1)}(\Mod_A)^{S^1} \\
&\simeq&  \lim_{\Spec C\in (\TSZ_A)_{/\widehat{S\times S}}}\Rep(\DD_\infty(C\otimes_AS^1))(\Mod_A)^{S^1} \\
&\simeq& \Rep_H^{\circlearrowleft\wedge}(\widehat{S\times_kS})^{S^1}  \\
&\simeq& \Rep(\mathbb{T}_{A/k}[-1]^{S^1})(\Mod_A^{S^1})
\end{eqnarray*}
where the second functor is induced by $\mathcal{I}:\QC\to \LMod\circ \DD_1$ (cf. Section~\ref{liemodule}), and the first equivalence and the second equivalence come from \cite[Proposition 7.1]{IMA} and \cite[Lemma 7.18]{DC1}.
According to Corollary~\ref{factorizationcomplete},
the second arrow factors as 
\begin{eqnarray*}
\lim_{\Spec C\in (\TSZ_A)_{/\widehat{S\times S}}}\QC(C\otimes_AS^1)^{S^1} &\to& \lim_{\Spec C\in (\TSZ_A)_{/\widehat{S\times S}}}\QC_H'((\Spec C\otimes_AS^1)^\wedge)^{S^1}=\QC_H^{\circlearrowleft\wedge}(\widehat{S\times_kS})^{S^1} \\
&\to&  \lim_{\Spec C\in (\TSZ_A)_{/\widehat{S\times S}}}\LMod_{\DD_1(C\otimes_AS^1)}(\Mod_A)^{S^1}.
\end{eqnarray*}
Thus, our claim follows.
\QED

For later use, we give the following definition.

\begin{Definition}
\label{lateruse}
For $B\in (\CCAlg_k)_{A//A}$
we define
\[
\textup{comp}_{B}^\circlearrowleft : \QC_H^\circlearrowleft(\widehat{\Spec B})\to \QC_H^{\circlearrowleft\wedge}(\widehat{\Spec B})
\]
to be the morphism in $\Fun(BS^1,\wCat)$ which is obtained
from $\textup{comp}_{\QC}\circ T:\QC\circ T  \to \QC_H'|_{(\CCAlg_k)_{A//A}} \circ T$
by passing to right Kan extensions along $\EXT_A\to (\CCAlg_k)_{A//A}$.
Taking $S^1$-invariants, we define the induced functor $(\textup{comp}_{B}^\circlearrowleft)^{S^1} : \QC_H^\circlearrowleft(\widehat{\Spec B})^{S^1}\to \QC_H^{\circlearrowleft\wedge}(\widehat{\Spec B})^{S^1}$.
\end{Definition}

\section{Revisiting Construction in Part I}

In this section, we revisit the construction of an
object of
\[
\QC_!(S)^{S^1}\times_{\QC_!((S\times_kS)_S^\wedge)^{S^1}}\QC_!((S\times_kLS)_S^\wedge)^{S^1}
\]
in \cite[Section 5,6,7]{DC1} in view of the results of this paper.

In {\it loc.cit.}, we denote the constructed object by $\mathcal{V}_\dagger$.
Consider the diagram
\[
\label{testdiagram}
\vcenter{\tag{7.1}
\xymatrix{
(S\times_kS)_S^\wedge \ar[r] \ar[d]^{\hat{\textup{pr}}_2} & (S\times_kLS)_S^\wedge \ar[d] \\
S \ar[r] & LS.
}
}
\]
The vertical morphisms are second projections.
The horizontal morphisms are determined by the morphism $S\to LS$ given by constant loops.
The $!$-pullback functors along morphisms in the diagram~\ref{testdiagram}
induces a categorical equivalence
\[
\QC_!(LS)^{S^1}\stackrel{\sim}{\to} \QC_!(S)^{S^1}\times_{\QC_!((S\times_kS)_S^\wedge)^{S^1}}\QC_!((S\times_kLS)_S^\wedge)^{S^1}.
\]
Through this equivalence, $\mathcal{V}_\dagger$  defines an object $\mathcal{H}_\circlearrowleft(\CCC)$ of $\QC_!(LS)^{S^1}$.

\vspace{2mm}

We will summarize the construction by using results of this paper
and highlighting several points relevant to the next section (cf. Theorem~\ref{main}).

(Step 1)
We first consider the image $\mathcal{V}'_\dagger$ of $\mathcal{V}_\dagger$ under the projection to $\QC_!((S\times_kLS)_S^\wedge)^{S^1}$.
We consider $\HH_\bullet(\CCC/A)$ to be the canonical $\mathbb{T}_S[-1]^{S^1}$-module, which is an object of 
$\Rep(\mathbb{T}_{S}[-1]^{S^1})(\Mod_A^{S^1})$ (see Section~\ref{modularsec}, \cite[Definition 6.3]{DC1}).

Let $\tilde{\mathcal{H}}$ be the object of $\lim_{\Spec C\in U_A}\QC(C\otimes_AS^1)^{S^1}\simeq \QC_H^{\circlearrowleft}(\widehat{S\times_kS})^{S^1}$,
defined in
Construction~\ref{Hconstruction} (see also Section~\ref{DE}).
By Proposition~\ref{modularmodule},
 $\HH_\bullet(\CCC/A)$ in $\Rep(\mathbb{T}_{S}[-1]^{S^1})(\Mod_A^{S^1})$ is naturally equivalent to the image of $\tilde{\mathcal{H}}$
under the composite 
\[
\QC_H^\circlearrowleft(\widehat{S\times_kS})^{S^1}\simeq \lim_{\Spec C\in U_A} \QC(C\otimes_AS^1)^{S^1}\stackrel{(\textup{comp}_{A\otimes A}^\circlearrowleft)^{S^1}}{\longrightarrow} \QC_H^{\circlearrowleft\wedge}(\widehat{S\times_kS})^{S^1}\stackrel{(\beta')^{S^1}}{\longrightarrow} \Rep(\mathbb{T}_{S}[-1]^{S^1})(\Mod_A^{S^1})
\]
(see Proposition~\ref{diagramprop}  and Definition~\ref{lateruse}).
In particular, $\HH_\bullet(\CCC/A)$ lies in the essential image of
the fully faithful functor
$(\beta')^{S^1}:\lim_{\Spec C\in U_A}\QC'_H(C\otimes_AS^1)^{S^1}\hookrightarrow \Rep(\mathbb{T}_{S}[-1]^{S^1})(\Mod_A^{S^1})$.
Let $\widehat{\mathcal{H}}$ be the image of $\tilde{\mathcal{H}}$ in $\QC_{H}^{\circlearrowleft\wedge}(\widehat{S\times_kS})^{S^1}=\lim_{\Spec C\in U_A}\QC'_H(C\otimes_AS^1)^{S^1}$, which can be identified with $\HH_\bullet(\CCC/A)$
endowed the canonical $\mathbb{T}_{S}[-1]^{S^1}$-action in $\Rep(\mathbb{T}_{S}[-1]^{S^1})(\Mod_A^{S^1})$.
Using the diagram
\[
\xymatrix{
\QC_!((S\times_kLS)^\wedge_S)^{S^1} & \QC_{H}^{\circlearrowleft\wedge}(\widehat{S\times_kS})^{S^1} \ar[l]_(0.4){\beta^{S^1}} \ar[r]^(0.45){(\beta')^{S^1}} & \Rep(\mathbb{T}_{S}[-1]^{S^1})(\Mod_A^{S^1})
}
\]
(see Proposition~\ref{diagramprop} (2)),
we define $\mathcal{V}'_\dagger$
to be the image of $\widehat{\mathcal{H}}$
in $\QC_!((S\times_kLS)^\wedge_S)^{S^1}$.

(Step 2)
By Proposition~\ref{diagramprop} (2), there exists the commutative diagram
\[
\label{diagramII}
\vcenter{
\xymatrix{
\lim_{\Spec C\in U_A}\QC_H(C\otimes_AS^1)^{S^1} \ar[r]^{\beta^{S^1}}  \ar[d] & \lim_{\Spec C\in U_A}\QC_!((\Spec C\otimes_AS^1)^\wedge)^{S^1}  \ar[d]  & \QC_!( (S\times_kLS)^\wedge_S)^{S^1} \ar[l]^(0.4){\simeq}   \ar[d] \\
\lim_{\Spec C\in U_A}\QC(C)^{S^1} \ar[r] & \lim_{\Spec C\in U_A}\QC_!(\Spec C)^{S^1} & \QC_!( (S\times_kS)^\wedge_S)^{S^1} \ar[l]^(0.45){\simeq}
}\tag{7.2}
}
\]
in $\wCat$.
Notice that there exists an $S^1$-equivariant canonical equivalence
\[
\Spec (A\otimes_kA)\otimes_AS^1\simeq \Spec A\times_k\Spec (A\otimes_kS^1)=S\times_kLS
\]
over $S=\Spec A$, where the structure morphism from the right side is the first projection.
Here the $A$-module structure of $A\otimes_k A$ in $(A\otimes_kA)\otimes_AS^1$
is defined by $A\simeq A\otimes_kk\to A\otimes_kA$ (that corresponds to the first projection
$S\times_kS\to S$). 
It follows that $(\Spec (A\otimes_kA)\otimes_AS^1)^\wedge\simeq (S\times_kLS)^\wedge_S$.
The upper equivalence in the diagram~\ref{diagramII} is the canonical functor
\[
\QC_!((\Spec (A\otimes_kA)\otimes_AS^1)^\wedge)^{S^1}\to \lim_{\Spec C\in U_A}\QC_!((\Spec C\otimes_AS^1)^\wedge)^{S^1}
\]
which is induced by the $!$-pullback functoriality over $\overline{U}_A$.
The lower equivalence is the canonical functor 
$\QC_!((\Spec A\otimes_kA)^\wedge)^{S^1}\to \lim_{\Spec C\in U_A}\QC_!((\Spec C)^\wedge)^{S^1}$
defined in the same way.
We conclude from the diagram~\ref{diagramII} that
the image of $\mathcal{V}_\dagger'$ in $\QC_!((S\times_kS)^\wedge_S)^{S^1}$
can be identified with the image of $\widehat{\mathcal{H}}$ under $\lim_{\Spec C\in U_A}\QC'_H(C\otimes_AS^1)^{S^1}\to \lim_{\Spec C\in U_A}\QC'_H(C)^{S^1}\simeq \lim_{\Spec C\in U_A}\QC(C)^{S^1}\to\QC_!((S\times_kS)^\wedge_S)^{S^1}$.

(Step 3)
We consider the image $\mathcal{V}''_\dagger$ of $\mathcal{V}'_\dagger$ in $\QC_!( (S\times_kS)^\wedge_S)^{S^1}$.
There exists an equivalence $\mathcal{V}_\dagger''\simeq \textup{pr}_2^!(\HH_\bullet(\CCC/A))$
in $\QC_!((S\times_kS)^\wedge_S)^{S^1}$ where $\textup{pr}_2^!:\QC_!(S)\to \QC_!(S\times_kS)\to  \QC_!((S\times_kS)^\wedge_S)$ is the $!$-pullback along the composite of the second projection
and the canonical morphism $(S\times_kS)^\wedge_S\to S\times_kS$ (see \cite[Lemma 6.7]{DC1}).
For later use (see the proof of Theorem~\ref{main}), we review the equivalence 
$\mathcal{V}_\dagger''\simeq \textup{pr}_2^!(\HH_\bullet(\CCC/A))$ in detail.

As observed in (Step 2), taking into account the commutative diagram~\ref{diagramII}, we see that
$\mathcal{V}''_\dagger$ is naturally equivalent to
the image of $\tilde{\mathcal{H}}$ under the composite
\[
\xi:\QC_H^\circlearrowleft(\widehat{S\times_kS})=\lim_{\Spec C\in U_A}\QC(C\otimes_AS^1)^{S^1}\to \QC_H'(\widehat{S\times_kS})^{S^1} \to \QC_!((S\times_kS)^\wedge_S)^{S^1}.
\]
From the diagram~\ref{diagramFF} and the natural transformation $\Upsilon:\QC\to \QC_!$
between functors $\CCAlgft_k\to \wCat$
we obtain the diagram
\[
\label{diagramJJ}
\vcenter{
\xymatrix{
\QC((A\otimes_kA)\otimes_AS^1)^{S^1}  \ar[r] \ar[d] & \QC(A\otimes_kA)^{S^1} \ar[r] \ar[d] &  \QC_!(A\otimes_kA)^{S^1} \ar[d] \\
\lim_{\Spec C\in U_A}\QC(C\otimes_AS^1)^{S^1} \ar[r] & \lim_{\Spec C\in U_A} \QC(C)^{S^1} \ar[r] & \lim_{\Spec C\in U_A}\QC_!(C)^{S^1}
}\tag{7.3}}
\]
in $\wCat$.
The composite of the lower functors is $\xi$ up to the lower right equivalence in the diagram~\ref{diagramII}.
The right vertical functor can be identified with the $!$-pullback functor
$\QC_!(A\otimes_kA)^{S^1}=\QC_!(S\otimes_kS)^{S^1}\to \QC_!((S\times_kS)^\wedge_S)^{S^1}$
along the canonical morphism $(S\times_kS)^\wedge_S\to S\times_kS$
up to the lower right equivalence in the diagram~\ref{diagramII}.
Recall that $\sigma_0:\overline{U}_A^{op} \to \Mod(\Mod_A^{S^1})_{\mathcal{H}}$
carries $A\otimes_kA$ to $\HH_\bullet(\textup{pr}_2^*(\CCC)/A)=\HH_\bullet(\Perf_A\otimes_k\CCC/A)$
with the equivalence $\HH_\bullet(\textup{pr}_2^*(\CCC)/A)\otimes_{(A\otimes_kA)\otimes_AS^1}A\simeq \HH_\bullet(\CCC/A)$ in $\Mod_A^{S^1}$ (see the proof of Proposition~\ref{cocarsection} for $\sigma_0$).
Recall that $\Perf_A^\otimes$-module structure on
$\textup{pr}_2^*(\CCC)$ in $\HH_\bullet(\textup{pr}_2^*(\CCC)/A)$
defined to be the restriction of the $\Perf_{A\otimes_kA}^\otimes$-module
$\textup{pr}_2^*(\CCC)$ along
$\Perf_A^\otimes\to \Perf_{A\otimes_kA}^\otimes$ given by $A\simeq A\otimes_kk\to A\otimes_kA$.
We will write $\underline{A}\otimes_kA$ (resp. $A\otimes_k\underline{A}$) for the $A$-module $A\otimes_kA$ given by
$\textup{id}\otimes(k\to A): A\otimes_kk\to A\otimes_kA$
(resp. $(k\to A)\otimes_k\textup{id}:k\otimes_kA\to A\otimes_kA$).
Remember that there exists a canonical equivalence between
 $\tilde{\mathcal{H}}$ and the image of $\HH_\bullet(\textup{pr}_2^*(\CCC)/A)$ in $\lim_{\Spec C\in U_A}\QC(C\otimes_AS^1)^{S^1}$ (cf. Corollary~\ref{semiorigin}).
Consequently, to obtain $\mathcal{V}_\dagger''\simeq \textup{pr}_2^!(\HH_\bullet(\CCC/A))$, it is enough to show that the image of  $\HH_\bullet(\textup{pr}_2^*(\CCC)/A)$
in $\QC_!(A\otimes_kA)^{S^1}$ is naturally equivalent to $\textup{pr}_2^!(\HH_\bullet(\CCC/A))$
where we think of $\HH_\bullet(\CCC/A)$ as the object of $\QC_!(S)^{S^1}$ defined by
the equivalence $\QC(S)^{S^1} \stackrel{\sim}{\to} \QC_!(S)^{S^1}$
(namely, we abuse notation by writing $\HH_\bullet(\CCC/A)$ for $\Upsilon_A(\HH_\bullet(\CCC/A))$).
Since $\textup{pr}_2^*(\HH_\bullet(\CCC/A))$ in $\QC(A\otimes_kA)^{S^1}$
maps to $\textup{pr}_2^!(\HH_\bullet(\CCC/A))$ in $\QC_!(A\otimes_kA)^{S^1}$,
it will suffice to construct an equivalence
$(A\otimes_kA)\otimes_{(\underline{A}\otimes_kA)\otimes_AS^1}\HH_\bullet(\textup{pr}_2^*(\CCC)/A) \simeq \textup{pr}_2^*(\HH_\bullet(\CCC/A))=(A\otimes_k\underline{A})\otimes_{A}\HH_\bullet(\CCC/A)\simeq A\otimes_k\HH_\bullet(\CCC/A)$
in $\QC(A\otimes_kA)^{S^1}$.
Since $(A\otimes_kA)\otimes_kS^1\to (\underline{A}\otimes_kA)\otimes_AS^1$
can be identified with $(A\otimes_kS^1)\otimes_k(A\otimes_kS^1)\to A\otimes_k(A\otimes_kS^1)$
induced by $A\otimes_kS^1\to A$ in the ``first'' term.
Similarly, $(\underline{A}\otimes_k A)\otimes_AS^1\to A\otimes_kA$
can be identified with $A\otimes_k(A\otimes_kS^1)\to A\otimes_kA$
induced by $A\otimes_kS^1\to A$ in the ``second'' term.
Thus, by Proposition~\ref{hochschildbasechange} (1)
\begin{eqnarray*}
\HH_\bullet(\textup{pr}_2^*(\CCC)/A) &\simeq&
\bigl((\underline{A}\otimes_kA)\otimes_AS^1\bigr)\otimes_{A\otimes_kS^1}\HH_\bullet(\CCC/k) \\ 
&\simeq& \bigl(A\otimes_k(A\otimes_kS^1)\bigr)\otimes_{A\otimes_kS^1} \HH_\bullet(\CCC/k)  \\
 &\simeq& A \otimes_k\HH_\bullet(\CCC/k)=p_{LS}^*(\HH_\bullet(\CCC/k)).
\end{eqnarray*}
It follows that
\[
\kappa:(A\otimes_kA)\otimes_{A\otimes_k(A\otimes_kS^1)}\HH_\bullet(\textup{pr}_2^*(\CCC)/A)\simeq
 (A\otimes_k\underline{A})\otimes_A(A\otimes_{(A\otimes_kS^1)}\HH_\bullet(\CCC/k)))  \simeq \textup{pr}_2^*(\HH_\bullet(\CCC/A))
\]
where the final equivalence comes from $A\otimes_{(A\otimes_kS^1)}\HH_\bullet(\CCC/k)\simeq\HH_\bullet(\CCC/A)$ (see Proposition~\ref{hochschildbasechange} (2)),
and
the first equivalence comes from base changes along the commutative diagram
in $\CAlg_k$
\[
\label{diagramLL}
\vcenter{
\xymatrix{
k\otimes_k(A\otimes_kS^1) \ar[r]_{(k\to A)\otimes\textup{id}} \ar[d]_{\textup{id}\otimes(A\otimes S^1\to  A)} & A\otimes_k(A\otimes_kS^1) \ar[d]^{\textup{id}\otimes(A\otimes S^1\to  A)}   \\
k\otimes_kA \ar[r]_{(k\to A)\otimes\textup{id}} &  A\otimes_kA,
}\tag{7.4}}
\]
and the equivalence $p_{LS}^*(\HH_\bullet(\CCC/k))\simeq \HH_\bullet(\textup{pr}_2^*(\CCC)/A)$.

(Step 4)
$\mathcal{V}_\dagger'\in \QC_!((S\times_kLS)_S^\wedge)^{S^1}$, $\HH_\bullet(\CCC/A)\in \QC_!(S)^{S^1}$,
and the equivalence $\mathcal{V}_\dagger''\simeq \textup{pr}_2^!(\HH_\bullet(\CCC/A))$ induced  by $\kappa$,
determine an object $\mathcal{V}_\dagger\in \QC_!(S)^{S^1}\times_{\QC_!((S\times_kS)_S^\wedge)^{S^1}}\QC_!((S\times_kLS)_S^\wedge)^{S^1}$.

\section{Comparison result}

We will prove the following comparison result.

\begin{Theorem}
\label{main}
The image of $\HH_\bullet(\CCC/k)\in \QC(LS)^{S^1}$ 
under
\[
\QC(LS)^{S^1}\stackrel{\Upsilon_{A\otimes_kS^1}}{\longrightarrow} \QC_!(LS)^{S^1}\simeq \QC_!(S)^{S^1}\times_{\QC_!((S\times_kS)_S^\wedge)^{S^1}}\QC_!((S\times_kLS)_S^\wedge)^{S^1}
\]
is equivalent to $\mathcal{V}_\dagger$. Namely, $\Upsilon_{A\otimes_kS^1}(\HH_\bullet(\CCC/k))\simeq \mathcal{H}_\circlearrowleft(\CCC)$.
\end{Theorem}

\Proof
Consider the commutative diagram
\[
\label{diagramKK}
\vcenter{
\xymatrix{
LS  & S\times_kLS \ar[l]_{p_{LS}} & (S\times_kLS)_S^\wedge \ar[l]  \\
S \ar[u]^\iota &  S\times_kS \ar[u]_{\textup{id}_S\times \iota} \ar[l]^{\textup{pr}_2} & (S\times_kS)^\wedge_S \ar[l]  \ar[u]
}\tag{8.1}}
\]
where $\textup{pr}_2$ and $p_{LS}$ are second projections,  and other horizontal morphisms
are canonical morphisms.
Note that the left square corresponds to the diagram~\ref{diagramLL}.
When combined with the diagrams~\ref{diagramII} and~\ref{diagramJJ},
the functors $\QC:\CCAlg_k\to \wCat$ and $\QC_!:\Fun(\CCAlgftec_k,\SSS)^{op}\to \wCat$
induce the commutative diagram in $\wCat$:
\[
\label{diagramMM}
\vcenter{
\xymatrix{
\QC(S)^{S^1} \ar[r] \ar[d]^{\textup{id}} & \QC(S\times_kS)^{S^1} \ar[d] & \QC(S\times_kLS)^{S^1} \ar[l] \ar[d] \\
\QC(S)^{S^1} \ar[r] \ar[d]^{\Upsilon_A} & \lim_{\Spec C\in U_A}\QC(\Spec C)^{S^1} \ar[d] & \lim_{\Spec C\in U_A}\QC(\Spec C\otimes_AS^1)^{S^1} \ar[l] \ar[d] \\
\QC_!(S)^{S^1} \ar[r] & \lim_{\Spec C\in U_A}\QC_!(\Spec C)^{S^1} & \lim_{\Spec C\in U_A}\QC_!((\Spec C\otimes_AS^1)^\wedge)^{S^1} \ar[l] \\
\QC_!(S)^{S^1} \ar[u]^{\textup{id}} \ar[r] & \QC_!((S\times_kS)^\wedge_S)^{S^1} \ar[u]^{\simeq} & \QC_!((S\times_kLS)^\wedge_S)^{S^1} \ar[l] \ar[u]^{\simeq}.
}\tag{8.2}}
\]
The composition of vertical functors (taking inverses of equivalences)
determines the commutative diagram
\[
\label{diagramNN}
\vcenter{
\xymatrix{
\QC(S)^{S^1} \ar[r]^(0.4){\textup{pr}_2^*} \ar[d]^{\Upsilon_A} & \QC(S\times_kS)^{S^1} \ar[d] & \QC(S\times_kLS)^{S^1} \ar[l]^{(\textup{id}\times\iota)^*} \ar[d] \\
\QC_!(S)^{S^1} \ar[r]^(0.4){\textup{pr}_2^!} & \QC_!((S\times_kS)^\wedge_S)^{S^1} & \QC_!((S\times_kLS)^\wedge_S)^{S^1} \ar[l]^{(\textup{id}\times\iota)^!},
}\tag{8.3}}
\]
which is naturally equivalent to the diagram obtained from the diagram~\ref{diagramKK} by $*$-pullback functors,
$!$-pullback functors, and $\Upsilon$.
The middle vertical functor is the composite $\QC(S\times_kS)^{S^1}\stackrel{\Upsilon_{A\otimes_kA}}{\to} \QC_!(S\times_kS)^{S^1}\to 
\QC_!((S\times_kS)^\wedge_S)$, where the second functor is the $!$-pullback functor
along the canonical morphism. The right vertical functor is defined in a similar way.
This diagram~\ref{diagramNN} induces
\[
\QC(LS)^{S^1} \to \QC(S)^{S^1}\times_{\QC(S\times_kS)^{S^1}}\QC(S\times_kLS)^{S^1}\to \QC_!(S)^{S^1}\times_{\QC_!((S\times_kS)^\wedge_S)^{S^1}}\QC_!((S\times_kLS)^\wedge_S)^{S^1}.
\]
Here, for ease of notation we write $FP$ and $FP_!$
for the fiber product in the middle and the fiber product on the right side, respectively. 
It will suffice to prove that there exists
an equivalence  between $\mathcal{V}_\dagger$ and the image of $\HH_\bullet(\CCC/k)\in \QC(LS)^{S^1}$ in $FP_!$.
 By definition, the image of $\HH_\bullet(\CCC/k)$ in $FP$
 is the data consisting of the pair $(A\otimes_{(A\otimes_kS^1)}\HH_\bullet(\CCC/k),A\otimes_k\HH(\CCC/k)) \in \QC(S)^{S^1}\times \QC(S\times_kLS)^{S^1}$
  together with the canonical equivalence
\[
 \omega:(A\otimes_kA)\otimes_{A\otimes_k(A\otimes_kS^1)}(A\otimes_k\HH(\CCC/k)) \simeq \textup{pr}_2^*(A\otimes_{(A\otimes_kS^1)}\HH_\bullet(\CCC/k))
\] 
in $\QC(S\times_kS)^{S^1}$, which is obtained from  the pullback functoriality over the diagram~\ref{diagramLL}.
Consider the object $V_\dagger$ of $FP$, defined by
the pair $(\HH_\bullet(\CCC/A), \HH_\bullet(\textup{pr}_2^*(\CCC)/A)) \in \QC(S)^{S^1}\times \QC(S\times_kLS)^{S^1}$
together with the equivalence
$\kappa:(A\otimes_kA)\otimes_{A\otimes_k(A\otimes_kS^1)}\HH_\bullet(\textup{pr}_2^*(\CCC)/A)\simeq \textup{pr}_2^*(\HH_\bullet(\CCC/A))$ in (Step 3) in the previous Section.
From the construction in (Step 3), the image of $V_\dagger$ in $FP_!$ is
naturally equivalent to $\mathcal{V}_\dagger$.
Therefore, it is enough to construct an equivalence betweeen
$V_\dagger$ and the image of $\HH_\bullet(\CCC/k)$ in $FP$.
To this end,
we will use 
equivalences
$f:A\otimes_{(A\otimes_kS^1)}\HH_\bullet(\CCC/k)\stackrel{\sim}{\to} \HH_\bullet(\CCC/A)$
(see Proposition~\ref{hochschildbasechange} (2))
and $g:p_{LS}^*(\HH(\CCC/k))\stackrel{\sim}{\to} \HH_\bullet(\textup{pr}_2^*(\CCC)/A)$
(see (Step 3)). It gives rise to the equivalence
\[
f\times g:(A\otimes_{(A\otimes_kS^1)}\HH_\bullet(\CCC/k),p_{LS}^*(\HH(\CCC/k))) \simeq (\HH_\bullet(\CCC/A), \HH_\bullet(\textup{pr}_2^*(\CCC)/A)).
\]
Unfolding the definition, the inverse of $\kappa$ is the composite of equivalences 
\begin{eqnarray*}
\textup{pr}_2^*(\HH_\bullet(\CCC/A)) &\stackrel{\textup{pr}_2^*f}{\longleftarrow}&   \textup{pr}_2^*(A\otimes_{(A\otimes_kS^1)}\HH_\bullet(\CCC/k))) \\
&\stackrel{\omega}{\longleftarrow}& (A\otimes_kA)\otimes_{A\otimes_k(A\otimes_kS^1)}(A\otimes_k\HH_\bullet(\CCC/k)) \\
&\stackrel{(\textup{id}\times \iota)^*g}{\longrightarrow}& (A\otimes_kA)\otimes_{A\otimes_k(A\otimes_kS^1)}\HH_\bullet(\textup{pr}_2^*(\CCC)/A).
\end{eqnarray*}
Thus, we have an equivalence $\kappa\circ (\textup{id}\times \iota)^*g\simeq \textup{pr}_2^*f\circ \omega$.
This equivalence and $f\times g$ define an equivalence between $V_\dagger$
and the image of $\HH_\bullet(\CCC/k)$ in $FP$.
\QED

Let $\Omega^\circ(\CCC)$ and $\Omega^{\bullet}(\CCC)$ be two $\ZZ/2\ZZ$-periodic right crystals (D-modules), which are constructed from $\HH_\bullet(\CCC/k)$ and $\mathcal{H}_{\circlearrowleft}(\CCC)$,
respectively, in \cite[Section 8]{DC1}. Both $\Omega^\circ(\CCC)$ and $\Omega^{\bullet}(\CCC)$
have the underlying $\ZZ/2\ZZ$-periodic complex $\mathcal{HP}_\bullet(\CCC/A)$, that is,
the periodic cyclic homology/complex.
By the equivalence $\Upsilon_{A\otimes S^1}(\HH_\bullet(\CCC/k))\simeq \mathcal{H}_{\circlearrowleft}(\CCC)$ in $\QC_!(LS)^{S^1}$ in Theorem~\ref{main}, we see:

\begin{Corollary}
There exists an equivalence $\Omega^\circ(\CCC)\simeq \Omega^{\bullet}(\CCC)$
of  $\ZZ/2\ZZ$-periodic right crystals (D-modules).
\end{Corollary}

\end{document}